\documentclass[12pt,oneside]{amsart}

\usepackage{amssymb}  
\usepackage{latexsym} 
\usepackage{comment}

\usepackage[all]{xy}
\usepackage{rotating}

\DeclareFontEncoding{OT2}{}{} 
\newcommand{\textcyr}[1]{%
 {\fontencoding{OT2}\fontfamily{cmr}\fontseries{m}\fontshape{n}\selectfont #1}}

\newcommand{\Sha}{{\mbox{\textcyr{Sh}}}}

\newcommand{\Z}{{\mathbb Z}}
\newcommand{\Q}{{\mathbb Q}}

\newcommand{\C}{{\mathbb C}}
\newcommand{\F}{{\mathbb F}}
\newcommand{\BP}{{\mathbb P}}
\newcommand{\BA}{{\mathbb A}}
\newcommand{\CO}{{\mathcal O}}

\newcommand{\CG}{{\mathcal G}}

\newcommand{\To}{\longrightarrow}

\newcommand{\tensor}{\otimes}

\newcommand{\GL}{\operatorname{GL}}

\newcommand{\Gal}{\operatorname{Gal}}
\newcommand{\Aut}{\operatorname{Aut}}

\newcommand{\Hom}{\operatorname{Hom}}

\newcommand{\Mor}{\operatorname{Mor}}
\newcommand{\Spec}{\operatorname{Spec}}
\newcommand{\Sel}{\operatorname{Sel}}
\newcommand{\Selh}{\mathop{\widehat{\operatorname{Sel}}}}
\newcommand{\Pic}{\operatorname{Pic}}
\newcommand{\Jac}{\operatorname{Jac}}
\newcommand{\Alb}{\operatorname{Alb}}
\newcommand{\NS}{\operatorname{NS}}

\newcommand{\Gm}{\mathbb{G}_m}

\newcommand{\isom}{\cong}
\newcommand{\inj}{\hookrightarrow}

\newcommand{\res}{\operatorname{res}}

\newcommand{\cov}{{\text{\rm f-cov}}}
\newcommand{\sol}{{\text{\rm f-sol}}}
\newcommand{\ab}{{\text{\rm f-ab}}}
\newcommand{\nab}{{\text{\rm $n$-ab}}}
\newcommand{\can}{{\text{\rm can}}}

\newcommand{\diw}{\operatorname{div}}

\newcommand{\tors}{{\text{\rm tors}}}

\newcommand{\im}{\operatorname{im}}

\newcommand{\Br}{\operatorname{Br}}
\newcommand{\inv}{\operatorname{inv}}

\newcommand{\surj}{\longrightarrow\!\!\!\!\!\!\!\!\longrightarrow}
\newcommand{\Ad}[2]{#1(\BA_{#2})_{\bullet}}
\newcommand{\Cov}{\mathop{{\mathcal C}ov}}
\newcommand{\Sol}{\mathop{{\mathcal S}ol}}
\newcommand{\Ab}{\mathop{{\mathcal A}b}}
\newcommand{\ev}{\operatorname{ev}}
\newcommand{\et}{{\text{\'et}}}

\newenvironment{Proof}{\par\noindent{\sc Proof:}}%
                      {\hspace*{\fill}\nobreak$\Box$\par}

   {\begin{list}{}{\settowidth{\labelwidth}{#1}%
                   \setlength{\leftmargin}{\labelwidth}%
                   \addtolength{\leftmargin}{\labelsep}}}%
   {\end{list}}

\newtheorem{Theorem}{Theorem}[section]
\newtheorem{Lemma}[Theorem]{Lemma}
\newtheorem{Proposition}[Theorem]{Proposition}
\newtheorem{Corollary}[Theorem]{Corollary}
\newtheorem{Conjecture}[Theorem]{Conjecture}

\theoremstyle{definition}
\newtheorem{Definition}[Theorem]{Definition}
\newtheorem{Examples}[Theorem]{Examples}
\newtheorem{Remark}[Theorem]{Remark}
\newtheorem{Question}[Theorem]{Question}

\numberwithin{equation}{section}

\begin{document}

\title[Finite descent and rational points on curves]%
      {Finite descent obstructions\\ and rational points on curves}

\author{Michael Stoll}
\address{School of Engineering and Science,
         Jacobs University Bremen
         P.O.Box 750561,
	 28725 Bremen, Germany.}
\email{m.stoll@jacobs-university.de}
\date{September 28, 2007}


\begin{abstract}
  Let $k$ be a number field and $X$ a smooth projective $k$-variety.
  In this paper, we study the information obtainable from descent via
  torsors under finite $k$-group schemes on the location of the $k$-rational
  points on~$X$ within the adelic points. 
  Our main result is that if a curve $C/k$ maps nontrivially into an
  abelian variety $A/k$ such that $A(k)$ is finite and $\Sha(k,A)$
  has no nontrivial divisible elements, then
  the information coming from finite abelian descent
  cuts out precisely the rational points of~$C$. 
  We conjecture that this is the case for all curves of genus~$\ge 2$.
  We relate finite descent
  obstructions to the Brauer-Manin obstruction; in particular, we prove
  that on curves, the Brauer set equals the set cut out by finite abelian
  descent. 
  Our conjecture therefore implies that the Brauer-Manin obstruction against
  rational points in the only one on curves.
\end{abstract}

\maketitle

\tableofcontents


\section{Introduction}

In this paper we explore what can be deduced about the set of rational
points on a curve (or a more general variety) from a knowledge of its finite 
\'etale coverings.

Given a smooth projective variety~$X$ over a number field~$k$ and a
finite \'etale, geometrically Galois covering $\pi : Y \to X$, standard descent
theory tells us that there are only finitely many twists $\pi_j : Y_j \to X$
of~$\pi$ such that $Y_j$ has points everywhere locally, and then
$X(k) = \coprod_j \pi_j(Y_j(k))$. Since $X(k)$ embeds into the adelic points
$X(\BA_k)$, we obtain restrictions on where the rational points on~$X$
can be located inside $X(\BA_k)$: we must have
\[ X(k) \subset \bigcup_j \pi_j(D_j(\BA_k)) =: X(\BA_k)^{\pi} \,. \]
Taking the information from all such finite \'etale coverings together,
we arrive at 
\[ X(\BA_k)^{\cov} = \bigcap_{\pi} X(\BA_k)^{\pi} \,. \]
Since the information we get cannot tell us more than on which connected
component a point lies at the infinite places, we make a slight modification
by replacing the $v$-adic component of $X(\BA_k)$ with its set of connected
components, for infinite places $v$. 
In this way, we obtain $\Ad{X}{k}$ and $\Ad{X}{k}^{\cov}$. 

We can be more restrictive in the kind of coverings we allow. We denote
the set cut out by restrictions coming from finite abelian coverings only by
$\Ad{X}{k}^{\ab}$ and the set cut out by solvable coverings by 
$\Ad{X}{k}^{\sol}$. Then we have the chain of inclusions
\[ X(k) \subset \overline{X(k)} \subset \Ad{X}{k}^{\cov}
        \subset \Ad{X}{k}^{\sol} \subset \Ad{X}{k}^{\ab} \subset \Ad{X}{k} \,,
\]
where $\overline{X(k)}$ is the topological closure of $X(k)$ in~$\Ad{X}{k}$,
see Section~\ref{CoCo} below.

It turns out that the set cut out by the information
coming from finite \'etale abelian coverings on a curve~$C$ coincides with the
`Brauer set', which is defined using the Brauer group of~$C$:
\[ \Ad{C}{k}^\ab = \Ad{C}{k}^{\Br} \]
This follows easily from the descent theory of Colliot-Th\'el\`ene and Sansuc;
see Section~\ref{BM}. It should be noted, however, that this result seems
to be new. It says that on curves, all the information coming from torsors
under groups of multiplicative type is already obtained from torsors under
finite abelian group schemes.

In this way, it becomes possible to
study the Brauer-Manin obstruction on curves via finite \'etale abelian
coverings. For example, we provide an alternative proof of the main
result in Scharaschkin's thesis~\cite{Scharaschkin} characterizing
$\Ad{C}{k}^{\Br}$ in terms of the topological closure of the Mordell-Weil
group in the adelic points of the Jacobian, see Cor.~\ref{Sch}

Let us call $X$ ``good'' if it satisfies $\overline{X(k)} = \Ad{X}{k}^{\cov}$
and ``very good'' if it satisfies $\overline{X(k)} = \Ad{X}{k}^{\ab}$.

Then another consequence is that the Brauer-Manin obstruction is the only
obstruction against rational points on a curve that is very good. 
More precisely, the Brauer-Manin
obstruction is the only one against a weak form of weak approximation, i.e.,
weak approximation with information at the infinite primes reduced to
connected components.

An abelian variety $A/k$ is very good if and only if the 
divisible subgroup of $\Sha(k, A)$ is trivial. For example, if $A/\Q$ is a
modular abelian variety of analytic rank zero, then $A(\Q)$ 
and $\Sha(\Q, A)$ are both finite, and $A$
is very good. A principal homogeneous space $X$ for~$A$ such that 
$X(k) = \emptyset$ is very good if and only if it represents a  non-divisible 
element of $\Sha(k, A)$. See Cor.~\ref{AV} and the text following~it.

The main result of this paper says that
if $C/k$ is a curve that has a nonconstant morphism $C \to X$,
where $X$ is (very) good and $X(k)$ is finite, then $C$ is (very) good
(and $C(k)$ is finite), see Prop.~\ref{Dom}. This implies that every 
curve~$C/\Q$ whose Jacobian has a nontrivial factor~$A$ that is a modular
abelian variety of analytic rank zero is very good, see Thm.~\ref{CorMor}.
As an application, we prove that all
modular curves $X_0(N)$, $X_1(N)$ and $X(N)$ (over~$\Q$) are very good,
see Cor.~\ref{ModCurves}. For curves without rational points, we have
the following corollary: \\[1mm]
{\em If $C/\Q$ has a non-constant morphism into a modular abelian variety
of analytic rank zero, and if $C(\Q) = \emptyset$, then the absence of
rational points is explained by the Brauer-Manin obstruction.}

This generalizes a result due to Siksek~\cite{Siksek} by removing all
assumptions related to the Galois action on the fibers of the morphism
over rational points.

\medskip

The paper is organized as follows. After a preliminary section setting up
notation, we prove in Section~\ref{AbVar} some results on abelian varieties, 
which will be needed later on, but are also interesting in themselves.
Then, in Section~\ref{TT}, we review torsors and twists and set up some
categories of torsors for later use. Section~\ref{CoCo} introduces the
sets cut out by finite descent information, as sketched above, and 
Section~\ref{CoCoRP} relates this to rational points. 
Next we study the relationship between our sets $\Ad{X}{k}^{\cov/\sol/\ab}$
and the Brauer set $\Ad{X}{k}^{\Br}$ and its variants. This is done in
Section~\ref{BM}. We then discuss
certain inheritance properties of the notion of being ``excellent''
(which is stronger than ``good'') in Section~\ref{PropEx}. This is then
the basis for the conjecture formulated and discussed in Section~\ref{Conjs}.

\subsection*{Acknowledgments}

I would like to thank Bjorn Poonen for very fruitful discussions
and Jean-Louis Colliot-Th\'el\`ene, Alexei Skorobogatov, David Harari
and the anonymous referee for reading previous versions
of this paper carefully and making some very useful comments
and suggestions. Further input was provided by Jordan Ellenberg, 
Dennis Eriksson,
Joost van Hamel, Florian Pop and Felipe Voloch. 
Last but not least,
thanks are due to the Centre \'Emile Borel of the Institut Henri Poincar\'e
in Paris for hosting a special semester on ``Explicit methods in number
theory'' in Fall~2004. A large part of the following has its origins
in discussions I had while I was there.


\section{Preliminaries}

In all of this paper, $k$ is a number field.

Let $X$ be a smooth projective variety over~$k$. We modify the definition 
of the set of adelic points of~$X$ in the following way.\footnote{This
notation was introduced by Bjorn Poonen.}
\[ \Ad{X}{k} = \prod_{v \nmid \infty} X(k_v)
                 \times \prod_{v \mid \infty} \pi_0(X(k_v)) \,.
\]
In other words, the factors at infinite places~$v$ are reduced
to the set of connected components of~$X(k_v)$. We then have a canonical
surjection $X(\BA_k) \surj \Ad{X}{k}$. Note that for a zero-dimensional
variety (or reduced finite scheme) $Z$, we have $Z(\BA_k) = \Ad{Z}{k}$.
We will occasionally be a bit sloppy
in our notation, pretending that canonical maps like $\Ad{X}{k} \to \Ad{X}{K}$
(for a finite extension $K \supset k$) or $\Ad{Y}{k} \to \Ad{X}{k}$ (for a
subvariety $Y \subset X$) are inclusions, even though they in general
are not at the infinite places. So for example, the intersection
$X(K) \cap \Ad{X}{k}$ means the intersection of the images of both sets
in $\Ad{X}{K}$.

If $X = A$ is an abelian variety over~$k$, then 
\[ \prod_{v \nmid \infty} \{0\} \times \prod_{v \mid \infty} A(k_v)^0 
     = A(\BA_k)_{\diw}
\]
is exactly the divisible subgroup of~$A(\BA_k)$. This implies that
\[ \Ad{A}{k}/n\Ad{A}{k} = A(\BA_k)/n A(\BA_k) \]
and then that
\[ \Ad{A}{k} = \lim\limits_{\longleftarrow} \Ad{A}{k}/n\Ad{A}{k}
             = \lim\limits_{\longleftarrow} A(\BA_k)/n A(\BA_k)
             = \widehat{A(\BA_k)}
\]
is (isomorphic to) its own component-wise pro-finite completion and also the 
component-wise pro-finite completion of the usual group of adelic points.

We will denote by $\widehat{A(k)} = A(k) \tensor_{\Z} \hat{\Z}$ the
pro-finite completion $\lim\limits_{\longleftarrow} A(k)/nA(k)$ of the
Mordell-Weil group~$A(k)$. By a result of Serre~\cite[Thm.~3]{Serre71}, 
the natural 
map $\widehat{A(k)} \to \widehat{A(\BA_k)} = \Ad{A}{k}$ 
is an injection and therefore
induces an isomorphism with the topological
closure $\overline{A(k)}$ of~$A(k)$ in~$\Ad{A}{k}$. 
We will re-prove this in Prop.~\ref{Inj} below,
and even show something stronger than that, see 
Thm.~\ref{InjMod}. (Our proof is based on a later result of Serre.)
Note that we have an exact sequence
\[ 0 \To A(k)_{\tors} \To \widehat{A(k)} \To \hat{\Z}^r \To 0 \,, \]
where $r$ is the Mordell-Weil rank of~$A(k)$; in particular,
\[ \widehat{A(k)}_{\tors} = A(k)_{\tors} \,. \]

Let $\Sel^{(n)}(k, A)$ denote the $n$-Selmer group of~$A$ over~$k$,
as usual sitting in an exact sequence
\[ 0 \To A(k)/n A(k) \To \Sel^{(n)}(k, A) \To \Sha(k, A)[n] \To 0 \,. \]
If $n \mid N$, we have a canonical map of exact sequences
\[ \SelectTips{cm}{}
   \xymatrix{ 0 \ar[r] & A(k)/N A(k) \ar[r] \ar[d] 
                       & \Sel^{(N)}(k, A) \ar[r] \ar[d]
                       & \Sha(k, A)[N] \ar[r] \ar[d]^{\cdot N/n} & 0 \\
              0 \ar[r] & A(k)/n A(k) \ar[r]
                       & \Sel^{(n)}(k, A) \ar[r]
                       & \Sha(k, A)[n] \ar[r] & 0
            }
\]
and we can form the projective limit
\[ \Selh(k, A) = \lim_{\longleftarrow} \Sel^{(n)}(k, A) \,, \]
which sits again in an exact sequence
\[ 0 \To \widehat{A(k)} \To \Selh(k, A) \To T\,\Sha(k, A) \To 0 \,, \]
where $T \Sha(k, A)$ is the Tate module of~$\Sha(k, A)$ (and exactness
on the right follows from the fact that the maps 
$A(k)/N A(k) \to A(k)/n A(k)$ are surjective). If $\Sha(k, A)$ is finite,
or more generally, if the divisible subgroup $\Sha(k, A)_{\text{div}}$
is trivial, then the Tate module vanishes, and $\Selh(k, A) = \widehat{A(k)}$.
Note also that since $T\,\Sha(k, A)$ is torsion-free, we have
\[ \Selh(k, A)_{\tors} = \widehat{A(k)}_{\tors} = A(k)_{\tors} \,. \]

By definition of the Selmer group, we get maps
\[ \Sel^{(n)}(k, A) \To A(\BA_k)/n A(\BA_k) = \Ad{A}{k}/n \Ad{A}{k} \]
that are compatible with the projective limit, so we obtain a canonical
map
\[ \Selh(k, A) \To \Ad{A}{k} \]
through which the map $\widehat{A(k)} \to \Ad{A}{k}$ factors.
We will denote elements of~$\Selh(k, A)$ by $\hat{P}$, $\hat{Q}$ and the like,
and we will write $P_v$, $Q_v$ etc.\ for their images in $A(k_v)$ or
$\pi_0(A(k_v))$, so that the map $\Selh(k, A) \to \Ad{A}{k}$ is written
$\hat{P} \longmapsto (P_v)_v$. (It will turn out that this map is
injective, see Prop.~\ref{Inj}.)

If $X$ is a $k$-variety, then we use notation like $\Pic_X$, $\NS_X$, etc.,
to denote the Picard group, N\'eron-Severi group, etc., of~$X$ over~$\bar{k}$,
as a $k$-Galois module.

Finally, we will denote the absolute Galois group of~$k$ by~$\CG_k$.


\section{Some results on abelian varieties} \label{AbVar}

In the following, $A$ is an abelian variety over~$k$ of dimension~$g$.
For $N \ge 1$, we set $k_N = k(A[N])$ for the $N$-division field, and
$k_\infty = \bigcup_N k_N$ for the division field.

The following lemma, based on a result due to Serre on the image of the
Galois group in $\Aut(A_{\tors})$, forms the basis for the results
of this section.

\begin{Lemma} \label{L3}
  There is some $m \ge 1$ such that $m$ kills all the cohomology groups
  $H^1(k_N/k, A[N])$.
\end{Lemma}

\begin{Proof}
  By a result of Serre~\cite[p.~60]{SerreLettre}, the
  image of $\CG_k$ in $\Aut(A_\tors) = \GL_{2g}(\hat{\Z})$ meets the
  scalars $\hat{\Z}^\times$ in a subgroup containing $S = (\hat{\Z}^\times)^d$
  for some $d \ge 1$. We can assume that $d$ is even.
  
  Now we note that in
  \[ H^1(k_N/k, A[N]) \stackrel{\text{inf}}{\inj} 
                 H^1(k_\infty/k, A[N]) \To
                 H^1(k_\infty/k, A_{\tors}) \,,
  \]
  the kernel of the second map is killed by $\#A(k)_{\tors}$. Hence it 
  suffices to show that $H^1(k_\infty/k, A_{\tors})$ is killed by some~$m$.
  
  Let $G = \Gal(k_\infty/k) \subset \GL_{2g}(\hat{\Z})$, 
  then $S \subset G$ is a normal subgroup.
  We have the inflation-restriction sequence
  \[ H^1(G/S, A_{\tors}^S) \To H^1(G, A_{\tors}) \To H^1(S, A_{\tors}) \,. \]
  Therefore it suffices to show that there is some integer $D \ge 1$ 
  killing both
  $A_{\tors}^S$ and $H^1(S, A_{\tors}) = H^1((\hat{\Z}^\times)^d, \Q/\Z)^{2g}$.
  
  For a prime $p$, we define
  \[ \nu_p = \min\{v_p(a^d-1) : a \in \Z_p^\times\} \,. \]
  It is easy to see that when $p$ is odd, we have $\nu_p = 0$ if
  $p-1$ does not divide~$d$, and $\nu_p = 1 + v_p(d)$ otherwise. Also,
  $\nu_2 = 1$ if $d$ is odd (which we excluded), and $\nu_2 = 2 + v_2(d)$
  otherwise. In particular,
  \[ D = \prod_p p^{\nu_p} \]
  is a well-defined positive integer.
  
  We first show that $A_{\tors}^S$ is killed by~$D$. We have
  \[
    A_{\tors}^S = \Bigl(\bigoplus_p (\Q_p/\Z_p)^{(\Z_p^\times)^d}\Bigr)^{2g}\,,
  \]
  and for an individual summand, we see that
  \begin{align*}
    (\Q_p/\Z_p)^{(\Z_p^\times)^d}
      &= \{ x \in \Q_p/\Z_p : (a^d-1)x = 0 \quad\forall a \in \Z_p^\times \} \\
      &= \{ x \in \Q_p/\Z_p : p^{\nu_p} x = 0 \}
  \end{align*}
  is killed by~$p^{\nu_p}\!$, whence the claim.
  
  Now we have to look at $H^1(S, A_{\tors})$. It suffices to consider
  $H^1((\hat{\Z}^\times)^d, \Q/\Z)$. We start with
  \[ H^1((\Z_p^\times)^d, \Q_p/\Z_p) = 0 \,. \]
  Too see this, note that $(\Z_p^\times)^d$ is pro-cyclic (for odd~$p$,
  $\Z_p^\times$ is already pro-cyclic, for $p = 2$, $\Z_2^\times$ is
  $\{\pm 1\}$ times a pro-cyclic group, and the first factor goes
  away under exponentiation by~$d$, since $d$ was assumed to be even);
  let $\alpha \in (\Z_p^\times)^d$ be a topological generator. By evaluating
  cocycles at~$\alpha$, we obtain an injection
  \[ H^1((\Z_p^\times)^d, \Q_p/\Z_p)
      \inj \frac{\Q_p/\Z_p}{(\alpha-1)(\Q_p/\Z_p)}
      = \frac{\Q_p/\Z_p}{p^{\nu_p}(\Q_p/\Z_p)} = 0 \,.
  \]
  We then can conclude that $H^1((\hat{\Z}^\times)^d, \Q_p/\Z_p)$ is
  killed by~$p^{\nu_p}$. To see this, write 
  \[ (\hat{\Z}^\times)^d = (\Z_p^\times)^d \times T \,, \]
  where $T = \prod_{q \neq p} (\Z_q^\times)^d$. Then, by inflation-restriction
  again, there is an exact sequence
  \[ 0 = H^1((\Z_p^\times)^d, \Q_p/\Z_p)
      \To H^1((\hat{\Z}^\times)^d, \Q_p/\Z_p)
      \To H^1(T, \Q_p/\Z_p)^{(\Z_p^\times)^d} \,, \]
  and we have (note that $T$ acts trivially on $\Q_p/\Z_p$)
  \[ H^1(T, \Q_p/\Z_p)^{(\Z_p^\times)^d}
      = \Hom(T, (\Q_p/\Z_p)^{(\Z_p^\times)^d}) \,.
  \]
  This group is killed by~$p^{\nu_p}\!$, since 
  $(\Q_p/\Z_p)^{(\Z_p^\times)^d}$ is. It follows that
  \[ H^1((\hat{\Z}^\times)^d, \Q/\Z)
      = \bigoplus_p H^1((\hat{\Z}^\times)^d, \Q_p/\Z_p)
  \]
  is killed by $D = \prod_p p^{\nu_p}$.
  
  We therefore find that $H^1(G, A_{\tors})$ is killed by~$D^2$, and
  that $H^1(k_N/k, A[N])$ is killed by~$D^2 \#A(k)_{\tors}$, for all~$N$.
\end{Proof}

\begin{Remark}
  A similar statement is proved for elliptic curves in~\cite[Prop.~7]{Viada}.
\end{Remark}

\begin{Lemma} \label{L4}
  For all positive integers~$N$, the map 
  \[ Sel^{(N)}(k, A) \To \Sel^{(N)}(k_N, A) \]
  has kernel killed by~$m$, where $m$ is the number from Lemma~\ref{L3}.
\end{Lemma}

\begin{Proof}
  We have the following commutative and exact diagram.
  \[ \SelectTips{cm}{}
     \xymatrix{           & 0 \ar[d] & 0 \ar[d] \\
                0 \ar[r]  & \ker \ar[d] \ar[r]
                                     & H^1(k_N/k, A[N]) \ar[d]^{\text{inf}} \\
                0 \ar[r] & \Sel^{(N)}(k, A) \ar[d] \ar[r] 
                                     & H^1(k, A[N]) \ar[d]^{\text{res}} \\
                0 \ar[r] & \Sel^{(N)}(k_N, A) \ar[r]
                                     & H^1(k_N, A[N])
              }
  \]
  So the kernel in question injects into
  $H^1(k_N/k, A[N])$, and by Lemma~\ref{L3}, this group is killed by~$m$. 
\end{Proof}

\begin{Lemma} \label{Dens}
  Let $Q \in \Sel^{(N)}(k, A)$, and let $n$ be the order of~$mQ$,
  where $m$ is the number from Lemma~\ref{L3}. Then the density of places
  $v$ of~$k$ such that $v$ splits completely in $k_N/k$ and such that
  the image of~$Q$ in~$A(k_v)/N A(k_v)$ is trivial is at most 
  $1/(n [k_N : k])$.
\end{Lemma}

\begin{Proof}
  By Lemma~\ref{L4}, the kernel of $\Sel^{(N)}(k, A) \to \Sel^{(N)}(k_N, A)$
  is killed by~$m$. Hence the order of the image of~$Q$ in~$\Sel^{(N)}(k_N, A)$
  is a multiple of~$n$, the order of~$mQ$.
  Now consider the following diagram for a place
  $v$ that splits in~$k_N$ and a place $w$ of~$k_N$ above it.
  {\small
  \[ \SelectTips{cm}{}
     \xymatrix{ \Sel^{(N)}(k, A) \ar[d] \ar[r]
                 & \Sel^{(N)}(k_N, A) \ar[d] \ar@{^(->}[r]
                 & H^1(k_N, A[N]) \ar[d] \ar@{=}[r]
                 & \Hom(G_{k_N}, A[N]) \ar[d] \\
                A(k_v)/N A(k_v) \ar[r]^-{\cong}
                 & A(k_{N,w})/N A(k_{N,w}) \ar@{^(->}[r]
                 & H^1(k_{N,w}, A[N]) \ar@{=}[r]
                 & \Hom(G_{k_{N,w}}, A[N])
              }
  \]
  }
  Let $\alpha$ be the image of~$Q$ in~$\Hom(G_{k_N}, A[N])$. Then the image
  of~$Q$ is trivial in $A(k_v)/N A(k_v)$ if and only if $\alpha$ restricts
  to the zero homomorphism on~$G_{k_{N,w}}$. This is equivalent to saying
  that $w$ splits completely in $L/k_N$, where $L$ is the fixed field
  of the kernel of~$\alpha$. Since the order of~$\alpha$ is a multiple
  of~$n$, we have $[L : k_N] \ge n$, and the claim now follows from
  the Chebotarev Density Theorem.
\end{Proof}

Recall the definition of $\Selh(k, A)$ and the natural maps 
\[ A(k) \inj \widehat{A(k)} \inj \Selh(k, A) \To \Ad{A}{k} \,, \]
where we denote the rightmost map by
\[ \hat{P} \longmapsto (P_v)_v \,. \]
Also recall that $\Selh(k, A)_{\tors} = A(k)_{\tors}$ under the
identification given by the inclusions above.

\begin{Lemma} \label{Sep}
  Let $\hat{Q}_1, \dots, \hat{Q}_s \in \Selh(k, A)$ be elements of
  infinite order, and let $n \ge 1$. 
  Then there is some~$N$ such that the images
  of $\hat{Q}_1, \dots, \hat{Q}_s$ in~$\Sel^{(N)}(k, A)$ all have
  order at least~$n$.
\end{Lemma}

\begin{Proof}
  For a fixed $1 \le j \le s$, consider $(n-1)! \hat{Q}_j \neq 0$.
  There is some $N_j$ such that the image of $(n-1)! \hat{Q}_j$
  in~$\Sel^{(N_j)}(k, A)$ is non-zero. This implies that the image
  of~$\hat{Q}_j$ has order at least~$n$. Because of the canonical
  maps $\Sel^{(lN_j)}(k, A) \to \Sel^{(N_j)}(k, A)$, this will also
  be true for all multiples of~$N_j$. Therefore, any $N$ that is
  a common multiple of all the $N_j$ will do.
\end{Proof}

\begin{Proposition} \label{ZeroDim}
  Let $Z \subset A$ be a finite subscheme of an abelian
  variety~$A$ over~$k$ such that $Z(k) = Z(\bar{k})$. 
  Let $\hat{P} \in \Selh(k, A)$
  be such that $P_v \in Z(k_v) = Z(k)$ for a set of places~$v$ of~$k$
  of density~$1$. Then $\hat{P}$ is in the image of~$Z(k)$ in~$\Selh(k, A)$.
\end{Proposition}

\begin{Proof} 
  We first show that $\hat{P} \in Z(k) + A(k)_{\tors}$. (Here and in
  the following, we identify $A(k)$ with its image in~$\Selh(k, A)$.)
  Assume the contrary. Then none of
  the differences $\hat{P} - Q$ for $Q \in Z(k)$ has finite order. 
  Let $n > \#Z(k)$, then by Lemma~\ref{Sep}, we can find a number~$N$ 
  such that the image of~$m(\hat{P}-Q)$ under 
  $\Selh(k, A) \to \Sel^{(N)}(k, A)$ 
  has order at least $n$, for all $Q \in Z(k)$. 
  
  By Lemma~\ref{Dens}, the density of places of~$k$ such that $v$ splits
  in $k_N$ and at least one of $\hat{P} - Q$ (for $Q \in Z(k)$) maps trivially
  into $A(k_v)/N A(k_v)$ is at most
  \[ \frac{\#Z(k)}{n [k_N : k]} < \frac{1}{[k_N : k]} \,. \]
  Therefore, there is a set of places~$v$ of~$k$ of positive density such
  that $v$ splits completely in~$k_N/k$ and such that none of $\hat{P} - Q$
  maps trivially into $A(k_v)/N A(k_v)$. This implies
  $P_v \neq Q$ for all $Q \in Z(k)$, contrary to the assumption on~$\hat{P}$
  and the fact that $Z(k_v) = Z(k)$.
    
  It therefore follows that $\hat{P} \in Z(k) + A(k)_{\tors} \subset A(k)$. 
  Take a finite place $v$ of~$k$ such that $P_v \in Z(k)$ (the set of such
  places has density~1 by assumption). Then $A(k)$ injects into~$A(k_v)$.
  But the image~$P_v$ of~$\hat{P}$ under $\Selh(k, A) \to A(k_v)$ is in~$Z$,
  therefore we must have $\hat{P} \in Z(k)$.
\end{Proof}

The following is a simple, but useful consequence.

\begin{Proposition} \label{Inj}
  If $S$ is a set of places of~$k$ of density~$1$, then
  \[ \Selh(k, A) \To \prod_{v \in S} A(k_v)/A(k_v)^0 \]
  is injective. (Note that $A(k_v)^0 = 0$ for $v$ finite.) In particular,
  \[ \widehat{A(k)} \To \prod_{v \in S} A(k_v)/A(k_v)^0 \]
  is injective, and the canonical map $\widehat{A(k)} \to \Ad{A}{k}$ induces
  an isomorphism between $\widehat{A(k)}$ and $\overline{A(k)}$, the 
  topological closure of~$A(k)$ in~$\Ad{A}{k}$.
\end{Proposition}

This is essentially Serre's result in~\cite[Thm.~3]{Serre71}.

\begin{Proof}
  Let $\hat{P}$ be in the kernel. Then we can apply Prop.~\ref{ZeroDim}
  with $Z = \{0\}$, and we find that $\hat{P} = 0$.
  
  In the last statement, it is clear that the image of the map
  is~$\overline{A(k)}$, whence the result.
\end{Proof}

From now on, we will identify $\Selh(k, A)$ with its image in~$\Ad{A}{k}$.
We then have a chain of inclusions
\[ A(k) \subset \overline{A(k)} \subset \Selh(k, A) \subset \Ad{A}{k}\,, \]
and
\[ \Selh(k, A)/\overline{A(k)} \cong T\Sha(k, A) \]
vanishes if and only if the divisible subgroup of~$\Sha(k, A)$ is trivial.

We can prove a stronger result than the above. For a finite place~$v$
of~$k$, we denote by $\F_v$ the residue class field at~$v$. If $v$
is a place of good reduction for~$A$, then it makes sense to speak
of~$A(\F_v)$, the group of $\F_v$-points of~$A$. There is a canonical map
\[ \Selh(k, A) \To A(k_v) \To A(\F_v) \,. \]

\begin{Lemma} \label{Mod}
  Let $0 \neq \hat{Q} \in \Selh(k, A)$. Then there is a
  set of (finite) places~$v$ of~$k$ (of good reduction for~$A$) of positive 
  density such that the image of~$\hat{Q}$ in~$A(\F_v)$ is non-trivial.
\end{Lemma}

\begin{Proof}
  First assume that $\hat{Q} \notin A(k)_{\tors}$. Then $m\hat{Q} \neq 0$, 
  so there is 
  some $N$ such that $m\hat{Q}$ has nontrivial image in~$\Sel^{(N)}(k, A)$ 
  (where $m$ is, as usual, the number from Lemma~\ref{L3}). 
  By Lemma~\ref{Dens}, we find that
  there is a set of places~$v$ of~$k$ of positive density such that
  $Q_v \notin N A(k_v)$. Excluding the finitely many places dividing
  $N \infty$ or of bad reduction for~$A$ does not change this density.
  For $v$ in this reduced set, we have 
  $A(k_v)/N A(k_v) \cong A(\F_v)/N A(\F_v)$, and so the image of~$\hat{Q}$
  in~$A(\F_v)$ is not in~$N A(\F_v)$, let alone zero.
  
  Now consider the case that $\hat{Q} \in A(k)_{\tors} \setminus \{0\}$. We know
  that for all but finitely many finite places~$v$ of good reduction,
  $A(k)_{\tors}$ injects into $A(\F_v)$, so in this case, the statement
  is even true for a set of places of density~$1$.
\end{Proof}

\begin{Remark}
  Note that the corresponding statement for points $Q \in A(k)$ is trivial;
  indeed, there are only finitely many finite places~$v$ of good reduction
  such that $Q$ maps trivially into~$A(\F_v)$. (Consider some projective
  model of~$A$; then $Q$ and $0$ are two distinct points in projective
  space. They will reduce to the same point mod~$v$ if and only if $v$
  divides certain nonzero numbers ($2 \times 2$ determinants formed with
  the coordinates of the two points).) The lemma above says that things
  can not go wrong too badly when we replace $A(k)$ by its completion
  $\widehat{A(k)}$ or even~$\Selh(k, A)$.
\end{Remark}

\begin{Theorem} \label{InjMod}
  Let $S$ be a set of finite places of~$k$ of good reduction for~$A$
  and of density~1. Then the canonical homomorphisms
  \[ \Selh(k, A) \To \prod_{v \in S} A(\F_v) \text{\quad and\quad}
     \widehat{A(k)} \To \prod_{v \in S} A(\F_v)
  \]
  are injective.
\end{Theorem}

\begin{Proof}
  Let $\hat{Q}$ be in the kernel. If $\hat{Q} \neq 0$, then by Lemma~\ref{Mod}, 
  there
  is a set of places~$v$ of positive density such that the image of~$\hat{Q}$
  in~$A(\F_v)$ is nonzero, contradicting the assumptions. So $\hat{Q} = 0$, and
  the map is injective.
\end{Proof}

For applications, it is useful to remove the condition in Prop.~\ref{ZeroDim}
that all points of~$Z$ have to be defined over~$k$.

\begin{Theorem} \label{ZeroDim1}
  Let $Z \subset A$ be a finite subscheme of an abelian
  variety~$A$ over~$k$. Let $\hat{P} \in \Selh(k, A)$
  be such that $P_v \in Z(k_v)$ for a set of places~$v$ of~$k$
  of density~$1$. Then $\hat{P}$ is in the image of~$Z(k)$ in~$\Selh(k, A)$.
\end{Theorem}

\begin{Proof}
  Let $K/k$ be a finite extension such that $Z(K) = Z(\bar{k})$.
  By Prop.~\ref{ZeroDim}, we have that the image of~$\hat{P}$
  in~$\Ad{A}{K}$ is in~$Z(K)$. This implies that the image of~$\hat{P}$
  in~$\Ad{A}{K}$ is in~$Z(k)$ (since $\hat{P}$ is $k$-rational). Now the
  canonical map $\Ad{A}{k} \to \Ad{A}{K}$ is injective except possibly at
  some of the infinite places, so $P_v \in Z(k)$ for all but finitely
  many places. Now, replacing $Z$ by $Z(k)$ and applying Prop.~\ref{ZeroDim}
  again (this time over~$k$), we find that $\hat{P} \in Z(k)$, as claimed.
\end{Proof}

We have seen that for zero-dimensional subvarieties $Z \subset A$,
we have $\Ad{Z}{k} \cap \overline{A(k)} = Z(k)$, or even more generally,
$\Ad{Z}{k} \cap \Selh(k, A) = Z(k)$ (writing intersections
for simplicity). One can ask if this is valid more generally for
subvarieties $X \subset A$ that do not contain the translate
of an abelian subvariety of positive dimension.

\begin{Question} \label{AML}
  Is there such a thing as an ``Adelic Mordell-Lang Conjecture''?
  
  A possible statement is as follows. Let $A/k$ be an abelian variety and
  $X \subset A$ a subvariety not containing the translate of a 
  nontrivial subabelian variety of~$A$. 
  Then there is a finite subscheme $Z \subset X$ such that
  \[ \Ad{X}{k} \cap \Selh(k,A) \subset \Ad{Z}{k} \,. \]
  
  If this holds, Thm.~\ref{ZeroDim1} above implies that 
  \[ X(k) \subset \Ad{X}{k} \cap \Selh(k,A)
          \subset \Ad{Z}{k} \cap \Selh(k,A) = Z(k) \subset X(k) 
  \]
  and therefore $X(k) = \Ad{X}{k} \cap \Selh(k,A)$. In the notation
  introduced in Section~\ref{CoCo} below and by the discussion in
  Section~\ref{CoCoRP}, this implies
  \[ X(k) \subset \Ad{X}{k}^\ab \subset \Ad{X}{k} \cap \Ad{A}{k}^\ab
          = \Ad{X}{k} \cap \Selh(k,A) = X(k) \,,
  \]
  and so $X$ is excellent w.r.t.\ abelian coverings (and hence ``very good'').
\end{Question}

\begin{Remark}
  Note that the Adelic Mordell Lang Conjecture formulated above is true
  when $k$ is a global function field, $A$ is ordinary, and $X$ is not
  defined over $k^p$ (where $p$ is the characteristic of~$k$), see
  Voloch's paper~\cite{VolochML}. (The result is also implicit 
  in~\cite{Hrushovski}.)  
\end{Remark}


\section{Torsors and twists} \label{TT}

In this section, we introduce the notions of torsors (under finite
\'etale group schemes) and twists, and we describe various constructions
that can be done with these objects.

Let $X$ be a smooth projective (reduced, but not necessarily geometrically 
connected) variety over~$k$.

We will consider the following category $\Cov(X)$. Its objects are
$X$-torsors $Y$
under~$G$ (see for example~\cite{Skorobogatov} for definitions), where $G$ is a
finite \'etale group scheme over~$k$. More concretely, the data
consists of a $k$-morphism $\mu : Y \times G \to Y$ describing a right action
of~$G$ on~$Y$, together with a finite \'etale
$k$-morphism $\pi : Y \to X$ such that the following diagram is cartesian
(i.e., identifies $Y \times G$ with the fiber product $Y \times_X Y$).
\[ \SelectTips{cm}{}
   \xymatrix{ Y \times G \ar[d]_{\text{pr}_1} \ar[r]^{\mu} & Y \ar[d]^{\pi} \\
              Y \ar[r]^{\pi} & X
            }
\]
We will usually just write $(Y, G)$ for such an object, with the maps
$\mu$ and~$\pi$ being understood.
Morphisms $(Y', G') \to (Y, G)$ in~$\Cov(X)$ are given by a pair of
maps ($k$-morphisms of (group) schemes) $\phi : Y' \to Y$ and 
$\gamma : G' \to G$ such that the obvious diagram commutes:
\[ \SelectTips{cm}{}
   \xymatrix{ Y' \times G' \ar[d]_{\phi \times \gamma} \ar[r]^-{\mu'}
                & Y' \ar[d]^\phi \ar[r]^{\pi'} & X \ar@{=}[d] \\
              Y \times G \ar[r]^-\mu & Y \ar[r]^\pi & X
            }
\]
Note that $\gamma$ is uniquely determined by~$\phi$: if $y' \in Y'$, 
$g' \in G'$, there is a unique $g \in G$ such that 
$\phi(y') \cdot g = \phi(y' \cdot g')$,
so we must have $\gamma(g') = g$.

We will denote by $\Sol(X)$ and $\Ab(X)$ the full subcategories of~$\Cov(X)$
whose objects are the torsors $(Y, G)$ such that $G$ is solvable or abelian,
respectively.

If $X' \to X$ is a $k$-morphism of (smooth projective) varieties, then
we can pull back $X$-torsors under~$G$ to obtain $X'$-torsors under~$G$.
This defines covariant functors $\Cov(X) \to \Cov(X')$, $\Sol(X) \to \Sol(X')$
and $\Ab(X) \to \Ab(X')$.

The following constructions are described for $\Cov(X)$, but they
are similarly valid for $\Sol(X)$ and $\Ab(X)$.

If $(Y_1, G_1), (Y_2, G_2) \in \Cov(X)$ are two $X$-torsors, then
we can construct their fiber product $(Y, G) \in \Cov(X)$, 
where $Y = Y_1 \times_X Y_2$ and $G = G_1 \times G_2$. 
More generally, if $(Y_1, G_1) \to (Y, G)$ and $(Y_2, G_2) \to (Y, G)$
are two morphisms in~$\Cov(X)$, there is a fiber product
$(Z, H) \in \Cov(X)$, where $Z = Y_1 \times_Y Y_2$ and
$H = G_1 \times_G G_2$.

If $(Y, G) \in \Cov(X)$ is an $X$-torsor, where now everything is over $K$ 
with a finite extension $K/k$, then we can apply restriction of scalars to
obtain $(R_{K/k} Y, R_{K/k} G) \in \Cov(R_{K/k} X)$.

If $(Y, G) \in \Cov(X)$ is an $X$-torsor and $\xi$ is
a cohomology class in $H^1(k, G)$, then we can construct the {\em twist}
$(Y_\xi, G_\xi)$ of~$(Y, G)$ by~$\xi$. Here $G_\xi$ is the inner form of~$G$ corresponding to~$\xi$ (compare, e.g., \cite[pp.~12, 20]{Skorobogatov}).
We will denote the structure maps by $\mu_\xi$ and~$\pi_\xi$.
Usually, $H^1(k, G)$ is just a pointed set with distinguished element
corresponding to the given torsor; if the torsor is abelian, $H^1(k, G)$ 
is a group, and $G_\xi = G$ for all $\xi \in H^1(k, G)$.

If $(\phi, \gamma) : (Y', G') \to (Y, G)$ is a morphism and $\xi \in H^1(k,G')$,
then we get an induced morphism 
$(Y'_\xi, G'_\xi) \to (Y_{\gamma_* \xi}, G_{\gamma_* \xi})$
(where $\gamma_*$ is the induced map $H^1(k, G') \to H^1(k, G)$).
Similarly, twists are compatible with pull-backs, fiber products
and restriction of scalars.

Twists are transitive in the following sense. If $(Y, G) \in \Cov(X)$ is an 
$X$-torsor and $\xi \in H^1(k, G)$, $\eta \in H^1(k, G_\xi)$, then
there is a $\zeta \in H^1(k, G)$ such that
$((Y_\xi)_\eta, (G_\xi)_\eta) \isom (Y_\zeta, G_\zeta)$. Conversely,
if $\xi$ and $\zeta$ are given, then there is an $\eta \in H^1(k, G_\xi)$
such that the relation above holds.

The following observation does not hold in general for $\Sol(X)$ and~$\Ab(X)$.
If $Y \stackrel{\pi}{\to} X$ is any finite \'etale morphism, then there
is some $(\tilde{Y}, G) \in \Cov(X)$ such that 
$\tilde{\pi} : \tilde{Y} \to X$ factors through~$\pi$. Also,
if we have $(Y, G) \in \Cov(X)$ and $(Z, H) \in \Cov(Y)$, then there
is some $(\tilde{Z}, \Gamma) \in \Cov(X)$ such that $\tilde{Z}$ maps
to~$Z$ over~$X$ and such that the induced map $\tilde{Z} \to Y$ gives
rise to a $Y$-torsor $(\tilde{Z}, \tilde{H}) \in \Cov(Y)$. This last
statement is also valid with $\Sol(X)$ and $\Sol(Y)$ in place of
$\Cov(X)$ and~$\Cov(Y)$ (since extensions of solvable groups are
solvable).


\section{Finite descent conditions} \label{CoCo}

In this section, we use torsors and their twists, as described in the
previous section, in order to obtain obstructions against rational points.
The use of torsors under finite abelian group schemes is classical;
it is what is behind the usual descent procedures on elliptic curves
or abelian varieties (and so one can claim that they go all the way back
to Fermat). The non-abelian case was first studied by
Harari and Skorobogatov~\cite{HarariSkorobogatov}; see also~\cite{Harari2000}.

The following theorem (going back to Chevalley and Weil \cite{ChevalleyWeil})
summarizes the standard facts about descent via torsors.
Compare also \cite[Lemma~4.1]{HarariSkorobogatov}
and~\cite[pp.~105,~106]{Skorobogatov}.

\begin{Theorem} \label{Descent}
  Let $(Y, G) \in \Cov(X)$ be a torsor, where $X$ is a smooth projective
  $k$-variety.
  \begin{enumerate}\addtolength{\itemsep}{2mm}
    \item $\displaystyle X(k) = \coprod_{\xi \in H^1(k,G)} \pi_\xi(Y_\xi(k))$.
    \item The {\em $(Y,G)$-Selmer set}
          \[ \Sel^{(Y,G)}(k, X)
               = \{\xi \in H^1(k,G) : \Ad{Y_\xi}{k} \neq \emptyset\}
          \]
          is finite: there are only finitely many twists $(Y_\xi, G_\xi)$
          such that $Y_\xi$ has points everywhere locally.
  \end{enumerate}
  At least in principle, the Selmer set in the second statement can be
  determined explicitly, and the union in the first statement can be 
  restricted to this finite set.
\end{Theorem}

The idea behind the following considerations is to see how much information
one can get out of the various torsors regarding the image of~$X(k)$
in~$\Ad{X}{k}$. Compare Definition~4.2 in~\cite{HarariSkorobogatov} and
Definition~5.3.1 in Skorobogatov's book~\cite{Skorobogatov}.

\begin{Definition}
  Let $(Y, G) \in \Cov(X)$ be an $X$-torsor.
  We say that a point $P \in \Ad{X}{k}$ {\em survives~$(Y, G)$}, if it lifts
  to a point in~$\Ad{Y_\xi}{k}$ for some twist $(Y_\xi, G_\xi)$ of~$(Y, G)$.
\end{Definition}

There is a cohomological description of this property. An $X$-torsor under~$G$
is given by an element of $H^1_{\text{\'et}}(X, G)$. Pull-back through the
map $\Spec k \to X$ corresponding to a point in~$X(k)$ gives a map
\[ X(k) \To H^1(k, G) \,. \]
Note that it is not necessary to refer to non-abelian \'etale cohomology here: 
the map $X(k) \to H^1(k, G)$ induced by a torsor $(Y, G)$ simply arises
by associating to a point $P \in X(k)$ its fiber $\pi^{-1}(P) \subset Y$,
which is a $k$-torsor under~$G$ and therefore corresponds to an element
of~$H^1(k, G)$.

We get a similar map on adelic points:
\[ \Ad{X}{k} \To \prod_v H^1(k_v, G) \]
There is the canonical restriction map
\[ H^1(k, G) \To \prod_v H^1(k_v, G) \,, \]
and the various maps piece together to give a commutative diagram:
\[ \SelectTips{cm}{}
   \xymatrix{ X(k) \ar[r] \ar[d] & H^1(k, G) \ar[d] \\
              \Ad{X}{k} \ar[r]   & \prod_v H^1(k_v, G)
            }
\]
A point $P \in \Ad{X}{k}$ survives $(Y, G)$ if and only if its image in
$\prod_v H^1(k_v, G)$ is in the image of the global set~$H^1(k, G)$.
The $(Y, G)$-Selmer set is then the preimage in~$H^1(k, G)$ of the image
of~$\Ad{X}{k}$; this is completely analogous to the definition of a Selmer 
group in case $X$ is an abelian variety~$A$, and $G = A[n]$ is the $n$-torsion
subgroup of~$A$.

Here are some basic properties.

\begin{Lemma} \label{PropSurv} \strut
  \begin{enumerate}\addtolength{\itemsep}{2mm}
    \item If $(\phi, \gamma) : (Y', G') \to (Y, G)$ is a morphism in~$\Cov(X)$, 
          and if $P \in \Ad{X}{k}$ survives $(Y', G')$, then $P$ also 
          survives~$(Y, G)$.
    \item If $(Y', G) \in \Cov(X')$ is the pull-back of $(Y, G) \in \Cov(X)$
          under a morphism $\psi : X' \to X$, then $P \in \Ad{X'}{k}$
          survives $(Y', G)$ if and only if $\psi(P)$ survives~$(Y, G)$.
    \item If $(Y_1, G_1), (Y_2, G_2) \in \Cov(X)$ have fiber product
          $(Y, G)$, then $P \in \Ad{X}{k}$ survives $(Y, G)$ if and only
          if $P$ survives both $(Y_1, G_1)$ and~$(Y_2, G_2)$.
    \item Let $X$ be over~$K$, where $K/k$ is a finite extension, and
          let $(Y, G) \in \Cov(X)$ be an $X$-torsor. Then
          $P \in \Ad{(R_{K/k} X)}{k}$ survives $(R_{K/k} Y, R_{K/k} G)$ if
          and only if its image in $\Ad{X}{K}$ survives~$(Y, G)$.
    \item If $(Y, G) \in \Cov(X)$ and $\xi \in H^1(k, G)$, then
          $P \in \Ad{X}{k}$ survives $(Y, G)$ if and only if $P$
          survives $(Y_\xi, G_\xi)$.
  \end{enumerate}
\end{Lemma}

\begin{Proof}
  \begin{enumerate}\addtolength{\itemsep}{2mm}
    \item By assumption, there are $\xi \in H^1(k, G')$ and 
          $Q \in \Ad{Y'_\xi}{k}$ such that $\pi'_\xi(Q) = P$. Now we have
          the morphism $\phi_\xi : Y'_\xi \to Y_{\gamma_* \xi}$ over~$X$, 
          hence $\pi_{\gamma_* \xi}(\phi_\xi(Q)) = \pi'_\xi(Q) = P$, whence
          $P$ survives~$(Y, G)$.
    \item Assume that $P$ survives~$(Y', G)$. Then there are 
          $\xi \in H^1(k, G)$ and
          $Q \in \Ad{Y'_\xi}{k}$ such that $\pi'_\xi(Q) = P$. There is
          a morphism $\Psi_\xi : Y'_\xi \to Y_\xi$ over~$\psi$, hence
          we have that $\pi_\xi(\Psi_\xi(Q)) = \psi(P)$, so $\psi(P)$
          survives~$(Y, G)$. Conversely, assume that $\psi(P)$ survives
          $(Y, G)$. Then there are $\xi \in H^1(k, G)$ and 
          $Q \in \Ad{Y_\xi}{k}$ such that $\pi_\xi(Q) = \psi(P)$. The
          twist $(Y'_\xi, G_\xi)$ is the pull-back of~$(Y_\xi, G_\xi)$
          under~$\psi$; in particular, $Y'_\xi = Y_\xi \times_X X'$, and
          so there is $Q' \in \Ad{Y'_\xi}{k}$ mapping to~$Q$ in~$Y_\xi$
          and to~$P$ in~$X'$. Hence $P$ survives~$(Y', G)$.
    \item We have obvious morphisms $(Y, G) \to (Y_i, G_i)$. So by part~(1),
          if $P$ survives $(Y, G)$, then it also survives $(Y_1, G_1)$
          and~$(Y_2, G_2)$.
          Now assume that $P$ survives both $(Y_1, G_1)$ and~$(Y_2, G_2)$.
          Then there are $\xi_1 \in H^1(k, G_1)$ and $\xi_2 \in H^1(k, G_2)$
          and points $Q_1 \in \Ad{Y_{1,\xi_1}}{k}$, 
          $Q_2 \in \Ad{Y_{2,\xi_2}}{k}$ such that $\pi_{1,\xi_1}(Q_1) = P$
          and $\pi_{2,\xi_2}(Q_2) = P$. Consider 
          $\xi = (\xi_1, \xi_2) \in H^1(k, G) = H^1(k, G_1) \times H^1(k, G_2)$.
          We have that $Y_\xi = Y_{1,\xi_1} \times_X Y_{2,\xi_2}$, hence
          there is $Q \in \Ad{Y_\xi}{k}$ mapping to $Q_1$ and~$Q_2$ under the
          canonical maps $Y_\xi \to Y_{i,\xi_i}$ ($i = 1,2$), and to~$P$
          under $\pi_\xi : Y_\xi \to X$. Hence $P$ survives~$(Y, G)$.
    \item We have $H^1(k, R_{K/k} G) = H^1(K, G)$, and the corresponding
          twists are compatible. For any $\xi$ in this set, we have
          $R_{K/k} Y_\xi = (R_{K/k} Y)_\xi$, and the adelic points
          $\Ad{(R_{K/k} Y_\xi)}{k}$ and $\Ad{Y_\xi}{K}$ are identified.
          The claim follows.
    \item This comes from the fact that every twist of $(Y, G)$ is also
          a twist of $(Y_\xi, G_\xi)$ and vice versa.
  \end{enumerate}
\end{Proof}

By the Descent Theorem~\ref{Descent}, it is clear that (the image 
in~$\Ad{X}{k}$ of) a rational point~$P \in X(k)$ survives every torsor.
Therefore it makes sense to study the set of adelic points that survive
every torsor (or a suitable subclass of torsors) in order to obtain 
information on the location of the rational points within the adelic points.
Note that the set of points in~$\Ad{X}{k}$ surviving a given torsor is
closed --- it is a finite union of images of compact sets $\Ad{Y_\xi}{k}$
under continuous maps.

We are led to the following definitions.

\begin{Definition} Let $X$ be a smooth projective variety over~$k$.
  \begin{enumerate}\addtolength{\itemsep}{2mm}
    \item $ \Ad{X}{k}^\cov =
                \{P \in \Ad{X}{k}
                   : \text{$P$ survives all $(Y,G) \in \Cov(X)$}\}\,.
          $
    
    \item $ \Ad{X}{k}^\sol =
               \{P \in \Ad{X}{k}
                  : \text{$P$ survives all $(Y,G) \in \Sol(X)$}\}\,. 
          $
          
    \item $ \Ad{X}{k}^\ab =
               \{P \in \Ad{X}{k}
                  : \text{$P$ survives all $(Y,G) \in \Ab(X)$}\}\,. 
          $
  \end{enumerate}
\end{Definition}
(The ``f'' in the superscripts stands for ``finite'', since we are dealing
with torsors under finite group schemes only.)

By the remark made before the definition above, we have
\[ X(k) \subset \overline{X(k)} \subset \Ad{X}{k}^\cov \subset \Ad{X}{k}^\sol
        \subset \Ad{X}{k}^\ab \subset \Ad{X}{k} \,.
\]
Here, $\overline{X(k)}$ is the topological closure of~$X(k)$ in~$\Ad{X}{k}$.

Recall the ``evaluation map'' for $P \in \Ad{X}{k}$ and $G$ a finite
\'etale $k$-group scheme,
\[ \ev_{P,G} : H^1_{\text{\'et}}(X, G) \To \prod_v H^1(k_v, G) \]
(the set on the left can be considered as the set of isomorphism classes
of $X$-torsors under~$G$) and the restriction map
\[ \res_G : H^1(k, G) \To \prod_v H^1(k_v, G) \,. \]
In these terms, we have
\[ \Ad{X}{k}^\cov
     = \bigcap_G \{P \in \Ad{X}{k} : \im(\ev_{P,G}) \subset \im(\res_G)\} \,,
\]
where $G$ runs through all finite \'etale $k$-group schemes. We obtain
$\Ad{X}{k}^\sol$ and $\Ad{X}{k}^\ab$ in a similar way, by restricting $G$
to solvable or abelian group schemes.

In the definition above, we can restrict to $(Y,G)$ with $Y$ connected
(over $k$) if $X$ is connected: if we have $(Y,G)$
with $Y$ not connected, then let $Y_0$ by a connected component of~$Y$,
and let $G_0 \subset G$ be the stabilizer of this component. Then
$(Y_0,G_0)$ is again a torsor of the same kind as~$(Y,G)$, and we have
a morphism $(Y_0,G_0) \to (Y,G)$. Hence, by Lemma~\ref{PropSurv},~(1),
if $P$ survives $(Y_0,G_0)$, then it also survives~$(Y,G)$.

However, we cannot restrict to geometrically connected torsors when
$X$ is geometrically connected. The reason is that there can be
obstructions coming from the fact that a suitable geometrically connected
torsor does not exist.

\begin{Lemma} \label{LConn}
  Assume that $X$ is geometrically connected. If there is a torsor
  $(Y, G) \in \Cov(X)$ such that $Y$ and all twists $Y_\xi$ are
  $k$-connected, but not geometrically connected, then 
  $\Ad{X}{k}^\cov = \emptyset$. The analogous statement holds for
  the solvable and abelian versions.
\end{Lemma}

\begin{Proof}
  If $Y_\xi$ is connected, but not geometrically connected, then
  $\Ad{Y_\xi}{k} = \emptyset$ (this is because the finite scheme
  $\pi_0(Y_\xi)$ is irreducible and therefore satisfies the Hasse
  Principle, compare the proof of Prop.~\ref{ZeroDim0}). Hence
  no point in~$\Ad{X}{k}$ survives~$(Y, G)$.
\end{Proof}

Let us briefly discuss how this relates to the geometric fundamental
group of~$X$ over~$\bar{k}$, assuming $X$ to be geometrically connected.
In the following, we write $\bar{X} = X \times_k \bar{k}$ etc.,
for the base-change of~$X$ to a variety over~$\bar{k}$.
Every torsor $(Y, G) \in \Cov(X)$ (resp., $\Sol(X)$ or $\Ab(X)$) gives rise
to a covering $\bar{Y} \to \bar{X}$ that is Galois with (solvable or abelian)
Galois group $G(\bar{k})$. The stabilizer~$\Gamma$ of a connected component
of~$\bar{Y}$ is then a finite quotient of the geometric fundamental group
$\pi_1(\bar{X})$. If we fix an embedding $k \to \C$, then $\pi_1(\bar{X})$
is the pro-finite completion of the topological fundamental group 
$\pi_1(X(\C))$, so $\Gamma$ is also a finite quotient of~$\pi_1(X(\C))$.
If $\Gamma$ is trivial, then $\pi_0(Y)$ is a $k$-torsor under~$G$,
and $(Y, G)$ is the pull-back of~$(\pi_0(Y), G)$ under the structure
morphism $X \to \Spec k$. We call such a torsor {\em trivial}. Note that
all points in $\Ad{X}{k}$ survive a trivial torsor (since their image
in $\Ad{(\Spec k)}{k} = (\Spec k)(k) = \{\text{pt}\}$ survives everything);
therefore trivial torsors do not give information. 

Conversely, given a finite quotient~$\Gamma$ of $\pi_1(\bar{X})$ or 
of~$\pi_1(X(\C))$, there is a corresponding covering $\bar{Y} \to \bar{X}$
that will be defined over some finite extension $K$ of~$k$. Let
$\pi : Y \to X_K$ be the covering over~$K$; it is a torsor under a
$K$-group scheme~$G$ such that $G(\bar{k}) = \Gamma$. We now construct
a torsor $(Z, R_{K/k} G) \in \Cov(X)$ that over~$K$ factors through~$\pi$.
By restriction of scalars, we obtain 
$(R_{K/k} Y, R_{K/k} G) \in \Cov(R_{K/k} X_K)$. We pull back via the
canonical morphism $X \to R_{K/k} X_K$ to obtain $(Z, R_{K/k} G) \in \Cov(X)$.
Over $K$, we have the following diagram.
\[ \SelectTips{cm}{}
   \xymatrix{ Z_K \ar[r] \ar[d]^{(R_{K/k} G)_K}
                & (R_{K/k} Y)_K \ar[r]^-{\can} \ar[d]^{(R_{K/k} G)_K}
                & Y \ar[d]^G \\
              X_K \ar[r]^-{\can} & (R_{K/k} X_K)_K \ar[r]^-{\can} & X_K
            }
\]
(Here the right hand horizontal maps come from the identity morphism
$W \to W$ of a $K$-variety~$W\!$, under the identification of
$\Mor_k(V, R_{K/k} W)$ with $\Mor_K(V_K, W)$, taking $V = R_{K/k} W$; for
$W = Y$ and $W = X_K$, respectively.)
The composition of the lower horizontal maps is the identity morphism,
hence $(Z_K, (R_{K/k} G)_K) \in \Cov(X_K)$ maps to $(Y, G)$. Note that
the torsor we construct is in $\Sol(X)$ (resp., $\Ab(X)$) when $\Gamma$ is
solvable (resp., abelian).

\begin{Lemma} \label{LemmaMapTwist}
  Let $X$ be geometrically connected, $(Y, G), (Y', G') \in \Cov(X)$
  such that $Y$ is geometrically connected and such that
  $(\bar{Y}, \bar{G})$ maps to $(\bar{Y}', \bar{G}')$ as torsors of~$\bar{X}$.
  Then there is a twist $(Y'_\xi, G'_\xi)$ of~$(Y', G')$ such that
  $(Y, G)$ maps to $(Y'_\xi, G'_\xi)$.
\end{Lemma}

\begin{proof}
  Let $(\phi, \gamma) : (\bar{Y}, \bar{G}) \to (\bar{Y}', \bar{G}')$ be the 
  given morphism. Note that by assumption, the covering maps
  $\pi : Y \to X$ and $\pi' : Y' \to X$ are defined over~$k$.
  For $\sigma \in \CG_k$, this implies that 
  $({}^\sigma \phi, {}^\sigma \gamma)$ is also a morphism
  $(\bar{Y}, \bar{G}) \to (\bar{Y}', \bar{G}')$. We can then
  consider the composite morphism
  \[ \bar{Y} \stackrel{(\phi, {}^\sigma \phi)}{\To}
       \bar{Y}' \times_{\bar{X}} \bar{Y}' \stackrel{\cong}{\To}
       \bar{Y}' \times \bar{G}' \stackrel{\operatorname{pr}_2}{\To} \bar{G'} \,.
  \]
  Since $\bar{Y}$ is connected and $\bar{G}'$ is discrete, this morphism
  must be constant; let $\xi_\sigma \in G'(\bar{k})$ be its image. 
  It can then be checked that $\xi = (\xi_\sigma)_{\sigma \in \CG_k}$
  is a $G'$-valued cocycle and that after twisting $(Y', G')$ by~$\xi$,
  the morphism $\phi$ becomes defined over~$k$; since $\gamma$ is
  uniquely determined by~$\phi$, the same is true for~$\gamma$.
\end{proof}

We still assume $X$ to be geometrically connected. Let us call a family
of torsors $(Y_i, G_i) \in \Cov(X)$ (resp.\ $\Sol(X)$ or $\Ab(X)$) 
with $Y_i$ geometrically connected a
{\em cofinal family of coverings of~$X$} (resp.\ of {\em solvable} or 
{\em abelian coverings of~$X$}) if for every (resp.\ every solvable or
abelian) connected $(\bar{Y}, \bar{G}) \in \Cov(\bar{X})$
(resp.\ $\Sol(\bar{X})$ or $\Ab(\bar{X})$), there is a torsor $(Y_i, G_i)$
such that $(\bar{Y}_i, \bar{G}_i)$ maps to $(\bar{Y}, \bar{G})$.
We then have the following.

\begin{Lemma} \label{LemmaCof}
  Let $X$ be geometrically connected.
  \begin{enumerate}\addtolength{\itemsep}{1mm}
    \item If $\Ad{X}{k}^\cov \neq \emptyset$, then there is a cofinal
          family of coverings of~$X$. A similar statement holds for
          $\Ad{X}{k}^\sol$ and solvable coverings, and for $\Ad{X}{k}^\ab$
          and abelian coverings.
    \item If $(Y_i, G_i)_i$ is a cofinal family of coverings of~$X$,
          then $P \in \Ad{X}{k}$ is in $\Ad{X}{k}^\cov$ if and only if
          $P$ survives every $(Y_i, G_i)$. Similarly for the solvable
          and abelian variants.
  \end{enumerate}
\end{Lemma}

\begin{Proof}
  \begin{enumerate}\addtolength{\itemsep}{2mm}
    \item Let $P \in \Ad{X}{k}^\cov$, and
          let $\bar{Y} \to \bar{X}$ be a finite \'etale Galois covering
          with Galois group~$\Gamma$.
          Then by the discussion before Lemma~\ref{LemmaMapTwist},
          there is a torsor $(Z, G) \in \Cov(X)$, which we can assume
          to be $k$-connected, such that
          $(\bar{Z}, \bar{G})$ maps to $(\bar{Y}, \Gamma)$.
          Without loss of generality (after perhaps twisting $(Z, G)$),
          we can assume that $(Z, G)$ lifts~$P$. This implies that
          $Z$ is geometrically connected (compare Lemma~\ref{LConn}).
          So if we take all torsors $(Z, G)$ obtained in this way, we
          obtain a cofinal family of coverings of~$X$. The proof in
          the solvable and abelian cases is analogous.
    \item The `only if' part is clear. So assume that $P$ survives
          all $(Y_i, G_i)$, and let $(Z, \Gamma) \in \Cov(X)$ be arbitrary.
          Let $\bar{Z}_0$ be a connected component of~$\bar{Z}$, and
          let $\bar{\Gamma}_0$ be the stabilizer of $\bar{Z}_0$. Then
          there is some $(Y_i, G_i)$ such that 
          $(\bar{Y}_i, \bar{G}_i) \to (\bar{Z}_0, \bar{\Gamma}_0)
             \to (\bar{Z}, \bar{\Gamma})$,
          hence by Lemma~\ref{LemmaMapTwist}, there is a twist 
          $(Z_\xi, \Gamma_\xi)$
          such that $(Y_i, G_i)$ maps to it. Since $P$ survives $(Y_i, G_i)$
          by assumption, it also survives $(Z_\xi, \Gamma_\xi)$ and therefore
          $(Z, \Gamma)$, by Lemma~\ref{PropSurv}. The proof in the solvable
          and abelian cases is again analogous.
  \end{enumerate}
\end{Proof}

\begin{Lemma} \label{LemmaPi1}
  Let $X$ be geometrically connected.
  \begin{enumerate}\addtolength{\itemsep}{1mm}
    \item If $\pi_1(\bar{X})$ is trivial (i.e., $X$ is simply connected),
          then $\Ad{X}{k}^\cov = \Ad{X}{k}$.
    \item If the abelianization $\pi_1(\bar{X})^{\text{\rm ab}}$ is
          trivial, then $\Ad{X}{k}^\ab = \Ad{X}{k}$.
    \item If $\pi_1(\bar{X})$ is abelian (resp., solvable),
          then $\Ad{X}{k}^\cov = \Ad{X}{k}^\ab$ (resp.,
          $\Ad{X}{k}^\cov = \Ad{X}{k}^\sol$).
  \end{enumerate}
\end{Lemma}

\begin{Proof}
  \begin{enumerate}\addtolength{\itemsep}{2mm}
    \item In this case, all torsors are trivial and are therefore survived
          by all points in~$\Ad{X}{k}$.
    \item Here the same holds for all abelian torsors.
    \item We always have $\Ad{X}{k}^\cov \subset \Ad{X}{k}^\ab$.
          So let $P \in \Ad{X}{k}^\ab$; then by Lemma~\ref{LemmaCof},~(1),
          there is a cofinal family $(Y_i, G_i)$ of abelian coverings
          of~$X$, and since $\pi_1(\bar{X})$ is abelian, this is also
          a cofinal family of coverings without restriction. By part~(2) of
          the same lemma, it suffices to check that $P$ survives
          all $(Y_i, G_i)$, which we know to be true, in order to conclude
          that $P \in \Ad{X}{k}^\cov$. Similarly for the solvable variant.
  \end{enumerate}
\end{Proof}

We now list some fairly elementary properties of the sets 
$\Ad{X}{k}^{\ab/\sol/\cov}$.

\begin{Proposition} \label{MorIncl}
  If $X' \stackrel{\psi}{\to} X$ is a morphism, then
  $\psi(\Ad{X'}{k}^\cov) \subset \Ad{X}{k}^\cov$. Similarly for the solvable
  and abelian variants.
\end{Proposition}

\begin{Proof}
  Let $P \in \Ad{X'}{k}^\cov$, and let $(Y, G) \in \Cov(X)$ be an $X$-torsor.
  By assumption, $P$ survives the pull-back $(Y', G)$ of~$(Y, G)$ under~$\psi$,
  so by Lemma~\ref{PropSurv}, part~(2), $\psi(P)$ survives~$(Y, G)$.
  Since $(Y, G)$ was arbitrary, $\psi(P) \in \Ad{X}{k}^\cov$. The same
  proof works for the solvable and abelian variants.
\end{Proof}

\begin{Lemma} \label{TwoPt}
  Let $Z = \Spec k \amalg \Spec k = \{P_1, P_2\}$. Then
  \[ \{P_1, P_2\} = Z(k) = \Ad{Z}{k}^\ab \,. \]
\end{Lemma}

\begin{Proof}
  Let $Q \in \Ad{Z}{k}$ and assume that $Q \notin Z(k)$. We have
  to show that $Q \notin \Ad{Z}{k}^\ab$. By assumption, there are
  places $v$ and $w$ of~$k$ such that $Q_v = P_1$ and $Q_w = P_2$.
  We will consider torsors under $G = \Z/2\Z$. Pick some $\alpha \in k^\times$
  such that $\alpha \notin (k_v^\times)^2$ and $\alpha \notin (k_w^\times)^2$.
  Let $Y = \Spec k(\sqrt{\alpha}) \amalg (\Spec k \amalg \Spec k)$;
  then $(Y, G) \in \Ab(Z)$ in an obvious way. We want to show that
  no twist $(Y_\xi, G)$ for $\xi \in H^1(k, G) = k^\times/(k^\times)^2$
  lifts~$Q$. Such a twist is of one of the following forms.
  \begin{align*}
    (Y_\xi, G) &= \Spec k(\sqrt{\alpha}) \amalg (\Spec k \amalg \Spec k) \\
    (Y_\xi, G) &= (\Spec k \amalg \Spec k) \amalg \Spec k(\sqrt{\alpha}) \\
    (Y_\xi, G) &= \Spec k(\sqrt{\beta}) \amalg \Spec k(\sqrt{\gamma})
  \end{align*}
  where in the last case, $\beta$ and $\gamma$ are independent
  in~$k^\times/(k^\times)^2$. In the first two cases, $Q$ does not lift,
  since in the first case, the first component does not lift~$Q_v$, and
  in the second case, the second component does not lift~$Q_w$ (by our
  choice of~$\alpha$). In the third case, there is a set of places of~$k$
  of density~$1/4$ that are inert in both $k(\sqrt{\beta})$ and
  $k(\sqrt{\gamma})$, so that $\Ad{Y_\xi}{k} = \emptyset$. In particular,
  $Q$ does not lift to any of these twists.
\end{Proof}

\begin{Proposition} \label{Disjoint}
  If $X = X_1 \amalg X_2 \amalg \dots \amalg X_n$ is a disjoint union,
  then 
  \[ \Ad{X}{k}^\cov = \coprod_{j=1}^n \Ad{X_j}{k}^\cov \,, \] 
  and similarly for the solvable and abelian variants.
\end{Proposition}

\begin{Proof}
  It is sufficient to consider the case $n = 2$. We have maps
  $X_1 \to X$ and $X_2 \to X$, so (by Prop.~\ref{MorIncl})
  $\Ad{X_1}{k}^\cov \amalg \Ad{X_2}{k}^\cov \subset \Ad{X}{k}^\cov$
  (same for $\cdot^\sol$ and $\cdot^\ab$). For the reverse inclusion,
  consider the morphism $X \to \Spec k \amalg \Spec k = Z$ mapping $X_1$
  to the first point and $X_2$ to the second point. If $Q \in \Ad{X}{k}^\ab$,
  then its image is in $\Ad{Z}{k}^\ab = Z(k)$ (by Prop.~\ref{MorIncl}
  again and Lemma~\ref{TwoPt}). This means that 
  $Q \in \Ad{X_1}{k} \amalg \Ad{X_2}{k}$. The claim then follows easily.
\end{Proof}

\begin{Proposition} \label{ZeroDim0}
  If $Z$ is a (reduced) finite scheme, then $\Ad{Z}{k}^\ab = Z(k)$.
\end{Proposition}

\begin{Proof}
  By Prop.~\ref{Disjoint}, it suffices to prove this when $Z = \Spec K$
  is connected. But in this case, it is known that $Z$ satisfies the
  Hasse Principle. On the other hand, if $Z(k) \neq \emptyset$, then
  $Z = \Spec k$ and $\Ad{Z}{k}$ has just one point, so $Z(k) = \Ad{Z}{k}$.
  
  (The statement that $\Spec K$ as a $k$-scheme satisfies the Hasse Principle
  comes down to the following fact: 
  
  {\em If a group~$G$ acts transitively on a finite set~$X$ such that every
  $g \in G$ fixes at least one element of~$X$, then $\#X = 1$.} 
  
  To see this, let $n = \#X$ and assume (w.l.o.g.) that $G \subset S_n$.
  The stabilizer~$G_x$ of $x \in X$ is a subgroup of index~$n$ in~$G$.
  By assumption, $G = \bigcup_{x \in X} G_x$, so 
  $G \setminus \{1\} = \bigcup_{x \in X} (G_x \setminus \{1\})$.
  Counting elements now gives $\#G - 1 \le n(\#G/n - 1) = \#G - n$,
  which implies $n = 1$.)
\end{Proof}

\begin{Remark}
  Note that the Hasse Principle does not hold in general for finite schemes.
  A typical counterexample is given by the $\Q$-scheme
  \[ \Spec \Q(\sqrt{13}) \amalg \Spec \Q(\sqrt{17})
       \amalg \Spec \Q(\sqrt{13 \cdot 17}) \,.
  \]
\end{Remark}

\begin{Proposition}
  We have 
  \[ \Ad{(X \times Y)}{k}^\cov = \Ad{X}{k}^\cov \times \Ad{Y}{k}^\cov \,. \]
  Similarly for the solvable and abelian variants.
\end{Proposition}
  
\begin{Proof}
  Prop.~\ref{MorIncl} implies that 
  \[ \Ad{(X \times Y)}{k}^\cov \subset \Ad{X}{k}^\cov \times \Ad{Y}{k}^\cov \]
  (and similarly for the solvable and abelian variants).
  
  For the other direction,   
  we can assume that $X$ and~$Y$ are $k$-connected, compare 
  Prop.~\ref{Disjoint}. If $X$ (say) is not geometrically connected, then
  $\Ad{X}{k} = \emptyset$, hence $\Ad{(X \times Y)}{k} = \emptyset$ as well,
  and the statement is trivially true. So we can assume that
  $X$ and~$Y$ are geometrically connected.
  
  We now use the fact that 
  $\pi_1(\bar{X} \times \bar{Y}) = \pi_1(\bar{X}) \times \pi_1(\bar{Y})$.
  Let $P \in \Ad{X}{k}^\cov$ and $Q \in \Ad{Y}{k}^\cov$.
  By Lemma~\ref{LemmaCof},~(1), there are cofinal families of coverings
  $(V_i, G_i)$ of~$X$ and $(W_j, H_j)$ of~$Y$, which we can assume to lift
  $P$, resp., $Q$. Then the products
  $(V_i \times W_j, G_i \times H_j)$ form a cofinal family of coverings
  of~$X \times Y$, and it is clear that they lift $(P, Q)$. By 
  Lemma~\ref{LemmaCof},~(2), this implies that 
  $(P, Q) \in \Ad{(X \times Y)}{k}^\cov$. 
  
  The solvable and abelian variants are proved similarly, using the
  corresponding product property of the maximal abelian and solvable 
  quotients of the geometric fundamental group.
\end{Proof}

\begin{Proposition} \label{Restrict}
  If $K/k$ is a finite extension and $X$ is a $K$-variety, then
  \[ \Ad{(R_{K/k} X)}{k}^\cov = \Ad{X}{K}^\cov \]
  (under the canonical identification $\Ad{(R_{K/k} X)}{k} = \Ad{X}{K}$),
  and similarly for the solvable and abelian variants.
\end{Proposition}

\begin{Proof}
  Let $P \in \Ad{(R_{K/k} X)}{k}^\cov$, and let $(Y, G) \in \Cov(X)$.
  By assumption, $P$ survives $(R_{K/k} Y, R_{K/k} G) \in \Cov(R_{K/k} X)$,
  so by Lemma~\ref{PropSurv}, part~(4), $P$ also survives~$(Y, G)$.
  Since $(Y, G)$ was arbitrary, $P \in \Ad{X}{K}^\cov$, so the left
  hand side is contained in the right hand side. 
  
  For the proof of the reverse inclusion, we can reduce to the
  case that $X$ is $K$-connected, by Prop.~\ref{Disjoint}.
  If $X$ is $K$-connected, but not geometrically connected, then
  $\Ad{(R_{K/k} X)}{k} = \Ad{X}{K} = \emptyset$, and there is nothing
  to prove. So we can assume that $X$ is geometrically connected.
  Take $P \in \Ad{X}{K}^\cov$. Then by Lemma~\ref{LemmaCof},
  there is a cofinal family $(Y_i, G_i)$ of coverings of~$X$. We show
  that $(R_{K/k} Y_i, R_{K/k} G_i)$ is then a cofinal family of coverings
  of~$R_{K/k} X$. Indeed, it is known that 
  $\overline{R_{K/k} X} \cong \bar{X}^{[K : k]}$ (with the factors
  coming from the various embeddings of~$K$ into~$\bar{k}$), so
  $\pi_1(\overline{R_{K/k} X}) \cong \pi_1(\bar{X})^{[K : k]}$.
  This easily implies the claim. Now, viewing $P$ as an element
  of~$\Ad{(R_{K/k} X)}{k}$, we see by Lemma~\ref{PropSurv} that
  $P$ survives every $(R_{K/k} Y_i, R_{K/k} G_i)$, hence
  $P \in \Ad{(R_{K/k} X)}{k}^\cov$.
  
  The same proof works for the solvable and abelian variants.
\end{Proof}

\begin{Proposition} \label{Extend}
  If $K/k$ is a finite extension, then 
  \[ \Ad{X}{k}^\cov \subset \Ad{X}{k} \cap \Ad{X}{K}^\cov \]
  and similarly for the solvable and abelian variants.
  Note that the intersection is to be interpreted as the pullback
  of $\Ad{X}{K}^\cov$ under the canonical map $\Ad{X}{k} \to \Ad{X}{K}$,
  which may not be injective at the infinite places.
\end{Proposition}

\begin{Proof}
  We have a morphism $X \to R_{K/k} X_K$, inducing the canonical map
  \[ \Ad{X}{k} \To \Ad{(R_{K/k} X_K)}{k} = \Ad{X}{K} \,. \]
  The claim now follows from combining Props. \ref{MorIncl} and~\ref{Restrict}.
\end{Proof}

We also have an analogue of the Descent Theorem~\ref{Descent}.

\begin{Proposition} \label{Cover}
  Let $(Y, G) \in \Cov(X)$ be an $X$-torsor. Then
  \[ \Ad{X}{k}^\cov = \bigcup \pi_\xi\bigl(\Ad{Y_\xi}{k}^\cov\bigr) \,, \]
  where the union is extended over all twists $(Y_\xi, G_\xi)$ of~$(Y, G)$,
  or equivalently, over the finite
  set of twists with points everywhere locally. A similar statement holds for
  the solvable variant, when $G$ is solvable.
\end{Proposition}

\begin{Proof}
  Note first that by Prop.~\ref{MorIncl}, the right hand side is a
  subset of the left hand side.
  
  For the reverse inclusion, take $P \in \Ad{X}{k}^\cov$. 
  To ease notation, we will suppress the group schemes when denoting
  torsors in the following.
  Let $Y_1, \dots, Y_s \in \Cov(X)$ (or $\Sol(X)$) 
  be the finitely many
  twists of~$Y$ such that $P$ lifts.  
  
  Define $\tau(j) \subset \{1, \dots, s\}$ to be the set of indices~$i$
  such that for every $X$-torsor $Z$ mapping to~$Y_j$ (or short: an
  $X$-torsor $Z$ over~$Y_j$), there is a twist
  $Z_\xi$ that lifts~$P$ and induces a twist of~$Y_j$ that is isomorphic
  to~$Y_i$. We make a number of claims about this function.
  
  (i) $\tau(j)$ is non-empty.
  
  To see this, note first that for any given~$Z$, the corresponding set 
  (call it $\tau(Z)$) is
  non-empty, since by assumption $P$ must lift to some twist of~$Z$,
  and this twist induces a twist of~$Y_j$ to which $P$ also lifts, hence
  this twist must be one of the~$Y_i$. Second, 
  if $Z$ maps to~$Z'$ (as $X$-torsors over~$Y_j$), 
  we have $\tau(Z) \subset \tau(Z')$. Third, for every pair of $X$-torsors
  $Z$ and~$Z'$ over~$Y_j$, their relative fiber product $Z \times_{Y_j} Z'$
  maps to both of them. Taking these together, we see that $\tau(j)$ is a
  filtered intersection of non-empty subsets of a finite set and hence
  non-empty.
  
  (ii) If $i \in \tau(j)$, then $\tau(i) \subset \tau(j)$.
  
  Let $h \in \tau(i)$, and let $Z$ be an $X$-torsor over~$Y_j$.
  By definition of~$\tau(j)$, there is a twist $Z_\xi$ of~$Z$ lifting~$P$
  and inducing
  the twist $Y_i$ of~$Y_j$. Now by definition of~$\tau(i)$, there is
  a twist $(Z_\xi)_\eta$ of~$Z_\xi$ lifting~$P$ and inducing the twist 
  $Y_h$ of~$Y_i$. By transitivity of twists, this means that we have 
  a twist of~$Z$ lifting~$P$ and inducing the twist $Y_h$ of~$Y_j$.
  Since $Z$ was arbitrary, this shows that $h \in \tau(j)$.
  
  (iii) For some~$j$, we have $j \in \tau(j)$.
  
  Indeed, selecting for each~$j$ some $\sigma(j) \in \tau(j)$ (this is
  possible by~(i)), the map $\sigma$ will have a cycle: $\sigma^m(j) = j$
  for some $m \ge 1$ and~$j$. Then by~(ii), it follows that $j \in \tau(j)$.
  
  For this specific value of~$j$, we have therefore proved that every
  $X$-torsor $Z$ over~$Y_j$ has a twist that lifts~$P$ and induces the
  trivial twist of~$Y_j$. This means in particular that this twist is
  also a twist of~$Z$ as a $Y_j$-torsor. 
  
  Now assume that $P$ does not
  lift to $\Ad{Y_j}{k}^\cov$ (or $\Ad{Y_j}{k}^\sol$). Since the preimages
  of~$P$ in~$\Ad{Y_j}{k}$
  form a compact set and since surviving a torsor is a closed condition,
  we can find a $Y_j$-torsor $V$ that is not survived by any of the
  preimages of~$P$. We can then find an $X$-torsor $Z$ mapping to~$V$,
  staying in $\Sol$ when working in that category. (Note that this step
  does not work for~$\Ab$, since extensions of abelian groups need not
  be abelian again.) But by what we have just proved, $Z$ has a twist
  as a $Y_j$-torsor that lifts a preimage of~$P$, a contradiction.
  Hence our assumption that $P$
  does not lift to $\Ad{Y_j}{k}^\cov$ (or $\Ad{Y_j}{k}^\sol$)
  must be false.
\end{Proof}

\begin{Remark}
  The analogous statement for $\Ad{X}{k}^\ab$ and $G$ abelian is not true
  in general: it would follow that $\Ad{X}{k}^\ab = \Ad{X}{k}^\sol$, but
  Skorobogatov (see~\cite[\S~8]{Skorobogatov} or~\cite{Skorobogatov99})
  has a celebrated example of a surface~$X$
  such that $\emptyset = \Ad{X}{k}^\sol \subsetneq \Ad{X}{k}^\ab$.
  In fact, there is an abelian covering $\pi : Y \to X$ such that 
  $\bigcup_\xi \pi_\xi(\Ad{Y_\xi}{k}^\ab) = \emptyset$, which therefore
  gives a counterexample to the abelian version of the statement.
  
  Skorobogatov shows that the ``Brauer set'' $\Ad{X}{k}^{\Br}$ is
  non\-empty. In a later paper~\cite[\S\,5.1]{HarariSkorobogatov}, Harari and 
  Skorobogatov
  show that there exists an obstruction coming from a nilpotent,
  non-abelian covering (arising from an abelian covering of~$Y$). 
  The latter means that $\Ad{X}{k}^\sol = \emptyset$,
  whereas the former implies that $\Ad{X}{k}^\ab \neq \emptyset$, since
  $\Ad{X}{k}^{\Br} \subset \Ad{X}{k}^\ab$; see Section~\ref{BM} below.
  The interest in this result comes from the fact that it is the first
  example known of a variety where there is no Brauer-Manin obstruction,
  yet there are no rational points.
\end{Remark}


\section{Finite descent conditions and rational points} \label{CoCoRP}

The ultimate goal behind considering the sets cut out in the
adelic points by the various covering conditions is to obtain
information on the rational points. There is a three-by-three
matrix of natural statements relating these sets, see the
diagram below. Here, $\overline{X(k)}$ is the topological closure
of~$X(k)$ in~$\Ad{X}{k}$.

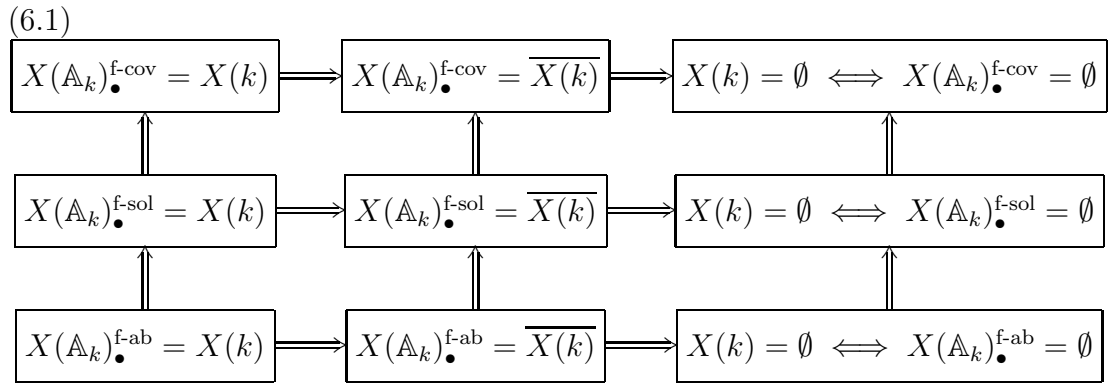
\begin{equation} \label{Props}
   \SelectTips{cm}{}
   \newcommand{\st}{\text{\Large\strut}}
   \xymatrix{ *+[F]{\Ad{X}{k}^\cov = X(k)\st} \ar@{=>}[r] &
                *+[F]{\Ad{X}{k}^\cov = \overline{X(k)}\st} \ar@{=>}[r] &
                *+[F]{X(k) = \emptyset \iff \Ad{X}{k}^\cov = \emptyset\st} \\
              *+[F]{\Ad{X}{k}^\sol = X(k)\st} \ar@{=>}[r] \ar@{=>}[u] &
                *+[F]{\Ad{X}{k}^\sol = \overline{X(k)}\st}
                    \ar@{=>}[r] \ar@{=>}[u] &
                *+[F]{X(k) = \emptyset \iff \Ad{X}{k}^\sol = \emptyset\st}
                 \ar@{=>}[u] \\
              *+[F]{\Ad{X}{k}^\ab = X(k)\st} \ar@{=>}[r] \ar@{=>}[u] &
                *+[F]{\Ad{X}{k}^\ab = \overline{X(k)}\st}
                    \ar@{=>}[r] \ar@{=>}[u] &
                *+[F]{X(k) = \emptyset \iff \Ad{X}{k}^\ab = \emptyset\st}
                  \ar@{=>}[u]
            }
\end{equation}

We have the indicated implications. If $X(k)$ is finite, then we
obviously have $X(k) = \overline{X(k)}$, and corresponding statements in
the left and middle columns are equivalent. In particular, this
is the case when $X$ is a curve of genus at least~2.

Let us discuss these statements. The ones in the middle column are
perhaps the most natural ones, whereas the ones in the left column
are better suited for proofs (as we will see below). The statements in the
right column can be considered as variants of the Hasse Principle;
in some sense they state that the Hasse Principle will eventually
hold if one allows oneself to replace $X$ by finite \'etale coverings.
Note that the weakest of the nine statements (the one in the upper
right corner), if valid for a class of varieties, would imply that
there is an effective procedure to decide whether there are $k$-rational
points on a variety~$X$ within that class or not: at least in principle,
we can list all the $X$-torsors and for each torsor
compute the finite set of twists with points everywhere locally. If this
set is empty, we know that $X(k) = \emptyset$. 
In order to obtain the torsors, we can for example enumerate all finite
extensions of the function field of~$X$ (assuming that $X$ is geometrically
connected, say) and check whether such an extension corresponds to an
\'etale covering of~$X$ that is a torsor under a finite group scheme.
On the other hand,
we can search for $k$-rational points on~$X$ at the same time, and
as soon as we find one such point, we know that $X(k) \neq \emptyset$.
The statement ``$X(k) = \emptyset \iff \Ad{X}{k}^\cov = \emptyset$''
guarantees that one of the two events must occur. (Note that
$\Ad{X}{k}^\cov$ can be written as a filtered intersection of compact
subsets of~$\Ad{X}{k}$, each coming from one specific torsor, so
if $\Ad{X}{k}^\cov = \emptyset$, then already one of these conditions
will provide an obstruction.)

For $X$ of dimension at least two, none of these statements can be
expected to hold in general. For example, a rational surface~$X$ has
trivial geometric fundamental group, and so $\Ad{X}{k}^\cov = \Ad{X}{k}$.
On the other hand, there are examples known of such surfaces that
violate the Hasse principle, so that we have 
$\emptyset = X(k) \subsetneq \Ad{X}{k}^\cov = \Ad{X}{k}$.
The first example (a smooth cubic surface) was given by
Swinnerton-Dyer~\cite{Swinnerton-Dyer}. There are also examples
among smooth diagonal cubic surfaces, see~\cite{CasselsGuy}, and
in~\cite{CT-Coray-Sansuc}, an infinite family of rational surfaces
violating the Hasse principle is given.

Let us give names to the properties in the left two columns in the
diagram~\ref{Props} above.

\begin{Definition}
  Let $X$ be a smooth projective $k$-variety. We call $X$
  \begin{enumerate}\addtolength{\itemsep}{1mm}
    \item {\em good with respect to all coverings} or simply {\em good} if
          $\overline{X(k)} = \Ad{X}{k}^\cov$,
    \item {\em good with respect to solvable coverings} if 
          $\overline{X(k)} = \Ad{X}{k}^\sol$,
    \item {\em good with respect to abelian coverings} or {\em very good} if
          $\overline{X(k)} = \Ad{X}{k}^\ab$,
    \item {\em excellent with respect to all coverings} if 
          $X(k) = \Ad{X}{k}^\cov$,
    \item {\em excellent with respect to solvable coverings} if 
          $X(k) = \Ad{X}{k}^\sol$,
    \item {\em excellent with respect to abelian coverings} if 
          $X(k) = \Ad{X}{k}^\ab$.
  \end{enumerate}
\end{Definition}


Now let us look at curves in more detail. When $C$ is a curve of genus~0,
then it satisfies the Hasse Principle, so
\[ \Ad{C}{k} = \emptyset \iff C(k) = \emptyset \,, \]
and then all the intermediate sets are equal and empty.
On the other hand, when $C(k) \neq \emptyset$, then $C \cong \BP^1$,
and $C(k)$ is dense in~$\Ad{C}{k}$, so
\[ \overline{C(k)} = \Ad{C}{k}^\cov = \Ad{C}{k}^\sol = \Ad{C}{k}^\ab = \Ad{C}{k} \,.
\]
So curves of genus~0 are always very good.

Now consider the case of a genus~1 curve. If $A$ is an elliptic
curve, or more generally, an abelian variety, then $\pi_1(\bar{A})$
is abelian, so by Lemma~\ref{LemmaPi1} we have
\[ \Ad{A}{k}^\cov = \Ad{A}{k}^\sol = \Ad{A}{k}^\ab \,. \]
Furthermore, among the abelian coverings, we can restrict to
the multiplication-by-$n$ maps $A \stackrel{n}{\to} A$. (In the terminology
used earlier, these coverings are a cofinal family.) This shows that
\[ \Ad{A}{k}^\ab = \Selh(k, A) \,. \]

Since the cokernel of the canonical map 
\[ \overline{A(k)} \cong \widehat{A(k)} \To \Selh(k, A) \]
is the Tate module of $\Sha(k, A)$, we get the following.

\begin{Corollary} \label{AV}
  \begin{align*}
    \text{$A$ is very good} &\iff \Sha(k, A)_{\diw} = 0 \\
    \text{$A$ is excellent w.r.t.\ abelian coverings}
       &\iff \text{$A(k)$ is finite and $\Sha(k, A)_{\diw} = 0$}
  \end{align*}
\end{Corollary}

See Wang's paper~\cite{Wang} for a discussion of the situation when one
works with $A(\BA_k)$ instead of~$\Ad{A}{k}$. Note that Wang's discussion
is in the context of the Brauer-Manin obstruction, which is closely related
to the ``finite abelian'' obstruction considered here, as discussed in 
Section~\ref{BM} below.

\begin{Corollary} \label{Ell0}
  If $A/\Q$ is a modular abelian variety of analytic rank zero, then 
  $A$ is excellent w.r.t.\ abelian coverings. 
  In particular, if $E/\Q$ is an elliptic
  curve of analytic rank zero, then $E$ is excellent w.r.t.\ abelian coverings.
\end{Corollary}

\begin{Proof}
  By work of Kolyvagin~\cite{Kolyvagin} and 
  Kolyvagin-Logachev~\cite{KolyvaginLogachev}, we know that $A(\Q)$ and 
  $\Sha(\Q, A)$ are both finite. By the above, it then follows that 
  $\Ad{A}{\Q}^\ab = A(\Q)$.
  
  If $E/\Q$ is an elliptic curve, then by work of Wiles~\cite{Wiles},
  Taylor-Wiles~\cite{TaylorWiles} and Breuil, Conrad, Diamond and 
  Taylor~\cite{BCDT}, we know that $E$ is modular and so the first assertion
  applies.
\end{Proof}

Now let $X$ be a principal homogeneous space for the abelian
variety~$A$. If $\Ad{X}{k} = \emptyset$, then all statements in~\eqref{Props}
are trivially true. So assume $\Ad{X}{k} \neq \emptyset$, and let
$\xi \in \Sha(k, A)$ denote the element corresponding to~$X$.
By Lemma~\ref{LemmaPi1}, we have
\[ \Ad{X}{k}^\cov = \Ad{X}{k}^\sol = \Ad{X}{k}^\ab \,, \]
and $\Ad{X}{k}^\ab = \emptyset$ if and only if $\xi \notin \Sha(k, A)_{\diw}$.
So for $\xi \neq 0$, $X$ is very good if and 
only if $\xi \notin \Sha(k, A)_{\diw}$ (since $X(k) = \emptyset$ in this case).

For curves~$C$ of genus~2 or higher, we always have that $C(k)$
is finite, and so the statements in the left and middle columns
in~\ref{Props} are equivalent. 
%
%
In this case, we can characterize the set $\Ad{C}{k}^\ab$ in a different way.

\begin{Theorem} \label{abGen}
  Let $C$ be a smooth projective geometrically connected curve over~$k$.
  Let $A = \Alb^0_C$ be its Albanese variety, and let $V = \Alb^1_C$ be
  the torsor under~$A$ that parametrizes classes of zero-cycles of degree~$1$
  on~$C$. Then there is a canonical map $\phi : C \to V$, and we have
  \[ \Ad{C}{k}^\ab = \phi^{-1}(\Ad{V}{k}^\ab) \,. \]
\end{Theorem}

Of course, since $C$ is a curve, $A$ is the same as the Jacobian variety
$\Jac_C = \Pic^0_C$, and $V$ is its torsor~$\Pic^1_C$, parametrizing
divisor classes of degree~$1$ on~$C$.

\begin{Proof}
  We know by Prop.~\ref{MorIncl} that 
  $\phi(\Ad{C}{k}^\ab) \subset \Ad{V}{k}^\ab$. It therefore suffices to
  prove that $\phi^{-1}(\Ad{V}{k}^\ab) \subset \Ad{C}{k}^\ab$.
  
  By~\cite[\S~VI.2]{SerreBook}, all (connected) finite abelian unramified
  coverings of~$\bar{C} = C \times_k \bar{k}$ are obtained through pull-back
  from isogenies into $\bar{V} \isom \bar{A}$. From this, we can deduce
  that the induced homomorphism 
  $\phi^* : H^1_\et(\bar{V}, \bar{G}) \to H^1_\et(\bar{C}, \bar{G})$
  is an isomorphism for all finite abelian $k$-group schemes~$G$. Since
  the map~$\phi$ is defined over~$k$, we obtain an isomorphism as $k$-Galois
  modules. The spectral sequence associated to the composition of functors
  $H^0(k, H^0_\et(\bar{V}, -)) = H^0_\et(V, -)$ (and
  similarly for $C$) gives a diagram with exact rows:
  \[ \SelectTips{cm}{}
     \xymatrix{ 0 \ar[r] & H^1(k, G) \ar[r] \ar@{=}[d] 
                         & H^1_\et(V, G) \ar[r] \ar[d]^{\phi^*}
                         & H^0(k, H^1_\et(\bar{V}, \bar{G})) \ar[r] 
                           \ar[d]_{\isom}^{\phi^*}
                         & H^2(k, G) \ar@{=}[d] \\
                0 \ar[r] & H^1(k, G) \ar[r]
                         & H^1_\et(C, G) \ar[r]
                         & H^0(k, H^1_\et(\bar{C}, \bar{G})) \ar[r]
                         & H^2(k, G)
              }
  \]
  By the 5-lemma, $\phi^* : H^1_\et(V, G) \to H^1_\et(C, G)$ is an isomorphism.
  
  Let $P \in \Ad{C}{k}$ such that $\phi(P) \in \Ad{V}{k}^\ab$, and let
  $(Y, G) \in \Ab(C)$. Then by the above, there is $(W, G) \in \Ab(V)$
  such that $Y$ is the pull-back of~$W$. By assumption, $\phi(P)$ survives
  $(W, G)$; without loss of generality, $(W, G)$ already lifts~$\phi(P)$.
  ($G$ is abelian, hence equal to all its inner forms.)
  Then $(Y, G)$ lifts~$P$, so $P$ survives~$(Y, G)$. Since $(Y, G)$ was
  arbitrary, $P \in \Ad{C}{k}^\ab$.
\end{Proof}

\begin{Remark} \label{RemAb}
  The result in the preceding theorem will hold more generally for
  smooth projective geometrically connected varieties~$X$ instead
  of curves~$C$, provided
  all finite \'etale abelian coverings of~$\bar{X}$ can be obtained
  as pullbacks of isogenies into the Albanese variety of~$X$. 
  For this, it is necessary and sufficient that the (geometric)
  N\'eron-Severi group of~$X$ is torsion-free, see~\cite[VI.20]{SerreBook}.
  
  For arbitrary varieties~$X$,
  we can define a set $\Ad{X}{k}^{\text{Alb}}$ consisting of the
  adelic points on~$X$ surviving all torsors that are pull-backs of
  $V$-torsors (where $V$ is the $k$-torsor under~$A$ that receives a canonical
  map~$\phi$ from~$X$), and then the result above will hold in the form
  \[ \Ad{X}{k}^{\text{Alb}} = \phi^{-1}(\Ad{V}{k}^\ab) \,. \]
  We trivially have $\Ad{X}{k}^\ab \subset \Ad{X}{k}^{\text{Alb}}$.

  In particular, we get that $\Ad{X}{k}^{\text{Alb}} = \Ad{X}{k}$ if $X$ 
  has trivial Albanese variety. For example, this is the case for all 
  complete intersections
  of dimension at least~$2$ in some projective space. (By Exercise~III.5.5
  in~\cite{Hartshorne}, $H^1(X, \CO) = 0$ in this case (over $\overline{k}$,
  say), so the Picard variety and therefore also its dual $\Alb^0(X)$ are
  trivial.) If in addition $\NS_X$ is torsion-free, then 
  $\Ad{X}{k}^\ab = \Ad{X}{k}$ as well.
\end{Remark}

\begin{Corollary} \label{abGen1}
  Let $C$ be a smooth projective geometrically connected curve over~$k$.
  Let $A$ be its Albanese (or Jacobian) variety, and let 
  $V = \Alb^1_C = \Pic^1_C$ as above.
  \begin{enumerate}\addtolength{\itemsep}{2mm}
    \item If $\Ad{C}{k} = \emptyset$, then $\Ad{C}{k}^\ab = C(k) = \emptyset$.
    \item If $\Ad{C}{k} \neq \emptyset$ and $V(k) \neq \emptyset$
          (i.e., $C$ has a $k$-rational divisor class of degree~$1$),
          then there is a $k$-defined embedding $\phi: C \inj A$, and we have
          \[ \Ad{C}{k}^\ab = \phi^{-1}(\Selh(k, A)) \,. \] 
          If $\Sha(k, A)_{\diw} = 0$, we have
          \[ \Ad{C}{k}^\ab = \phi^{-1}(\overline{A(k)}) \,. \]
    \item If $\Ad{C}{k} \neq \emptyset$ and $V(k) = \emptyset$,
          then, using the canonical map $\phi : C \to V$, we have
          \[ \Ad{C}{k}^\ab = \phi^{-1}(\Ad{V}{k}^\ab) \,. \]
          Let $\xi \in \Sha(k, A)$ be the element corresponding to~$V$.
          By assumption, $\xi \neq 0$. Then if $\xi \notin \Sha(k, A)_{\diw}$
          (and so in particular when $\Sha(k, A)_{\diw} = 0$),
          we have $C(k) = \Ad{C}{k}^\ab = \emptyset$.
  \end{enumerate}
  Similar statements are true for more general~$X$ in place of~$C$, with
  $\Ad{X}{k}^{\Alb}$ in place of~$\Ad{C}{k}^\ab$.
\end{Corollary}

\begin{Proof}
  This follows immediately from Thm.~\ref{abGen}, taking into account
  the descriptions of $\Ad{A}{k}^\ab$ and~$\Ad{V}{k}^\ab$ in Cor.~\ref{AV}
  and the text following it.
\end{Proof}

\bigskip

Let $X$ be a smooth projective geometrically connected $k$-variety,
let $A$ be its Albanese variety, and denote by $V$ the $k$-torsor under~$A$
such that there is a canonical map $\phi : X \to V$. ($V$ corresponds to
the cocycle class of $\sigma \mapsto [P^\sigma - P] \in A(\bar{k})$ for
any point $P \in X(\bar{k})$.) If $V(k) \neq \emptyset$, then $V$ is the trivial
torsor, and there is an $n$-covering of~$V$, i.e., a $V$-torsor under~$A[n]$.
So the non-existence of an $n$-covering of~$V$ is an obstruction against
rational points on~$V$ and therefore on~$X$.

If an $n$-covering of~$V$ exists, we can pull it back to a torsor
$(Y, A[n]) \in \Ab(X)$, and we will say that a point $P \in \Ad{X}{k}$
{\em survives the $n$-covering of~$X$} if it survives $(Y, A[n])$.
If there is no $n$-covering, then by definition no point in~$\Ad{X}{k}$
survives the $n$-covering of~$X$. If we denote the set of adelic points
surviving the $n$-covering of~$X$ by $\Ad{X}{k}^{\nab}$, then
we have
\[ \Ad{X}{k}^{\Alb} = \bigcap_{n \ge 1} \Ad{X}{k}^{\nab} \,. \]
In particular, for a curve $C$, we get
\[ \Ad{C}{k}^\ab  = \bigcap_{n \ge 1} \Ad{C}{k}^{\nab} \,. \]


\section{Relation with the Brauer-Manin obstruction} \label{BM}

In this section, we study the relationship between the finite covering
obstructions introduced in Section~\ref{CoCo} and the Brauer-Manin obstruction.
This latter obstruction was introduced by Manin~\cite{Manin} in~1970 in order
to provide a unified framework to explain violations of the Hasse Principle.

The idea is as follows. Let $X$ be, as usual, a smooth projective 
geometrically connected $k$-variety. We then have the (cohomological)
Brauer group
\[ \Br(X) = H^2_{\text{\'et}}(X, \Gm) \,. \]
If $K/k$ is any field extension and $P \in X(K)$ is a $K$-point of~$X$,
then the corresponding morphism $\Spec K \to X$ induces a homomorphism
$\phi_P : \Br(X) \to \Br(K)$. If $K = k_v$ is a completion of~$k$,
then there is a canonical injective homomorphism
\[ \inv_v : \Br(k_v) \inj \Q/\Z \]
(which is an isomorphism when $v$ is a finite place). In this way,
we can set up a pairing
\[ \Ad{X}{k} \times \Br(X) \To \Q/\Z\,, \quad
   ((P_v), b) \longmapsto
      \langle (P_v), b \rangle_{Br} = \sum_v \inv_v \bigl(\phi_{P_v}(b)\bigr) 
      \,.
\]
By a fundamental result of Class Field Theory, $k$-rational points on~$X$
pair trivially with all elements of~$\Br(X)$. This implies that
\[ \overline{X(k)} \subset \Ad{X}{k}^{\Br}
      = \{P \in \Ad{X}{k} :
            \langle P, b \rangle_{\Br} = 0 \text{\ for all $b \in \Br(X)$}\} \,.
\]
The set $\Ad{X}{k}^{\Br}$ is called the {\em Brauer set} of~$X$. If it is
empty, one says that there is a {\em Brauer-Manin obstruction} against
rational points on~$X$. More generally, if $B \subset \Br(X)$ is a
subgroup (or subset), we can define $\Ad{X}{k}^B$ in a similar way as
the subset of points in~$\Ad{X}{k}$ that pair trivially with all $b \in B$.

\medskip

The main result of this section is that for a curve~$C$, we have
\[ \Ad{C}{k}^{\Br} = \Ad{C}{k}^{\ab} \,, \]
see Cor.~\ref{BMabC} below. This implies that all the results we have deduced
or will deduce about finite abelian descent obstructions on curves also apply 
to the Brauer-Manin obstruction.

\bigskip

We first recall that the (algebraic) Brauer-Manin obstruction is
at least as strong as the obstruction coming from finite abelian descent.
For a more precise statement, see~\cite[Thm.~4.9]{HarariSkorobogatov}.
We define
\[ \Br_1(X) = \ker\bigl(\Br(X) \To \Br(X \times_k \bar{k})\bigr)
       \subset \Br(X)
\]
and set $\Ad{X}{k}^{\Br_1} = \Ad{X}{k}^{\Br_1(X)}$.

\begin{Theorem} \label{BMincl}
  For any smooth projective geometrically connected variety~$X$, we have 
  \[ \Ad{X}{k}^{\Br} \subset \Ad{X}{k}^{\Br_1} \subset \Ad{X}{k}^\ab \,. \]
\end{Theorem}
\begin{Proof}
  The main theorem of descent theory of Colliot-Th\'el\`ene and 
  Sansuc~\cite{CT-Sansuc}, as extended by Skorobogatov 
  (see \cite{Skorobogatov99} and~\cite[Thm.~6.1.1]{Skorobogatov}), states that
  $\Ad{X}{k}^{\Br_1}$ is equal to the set obtained from descent obstructions
  with respect to torsors under $k$-groups $G$ of multiplicative type,
  which includes all finite abelian $k$-groups. This proves the second
  inclusion. The first one follows from the definitions.
\end{Proof}

It is known that (see~\cite[Cor.~2.3.9]{Skorobogatov}; use that
$H^3(k, \bar{k}^\times) = 0$)
\[ \frac{\Br_1(X)}{\Br_0(X)} \isom H^1(k, \Pic_X) \,, \]
where $\Br_0(X)$ denotes the image of $\Br(k)$ in~$\Br(X)$.
We also have the canonical map $H^1(k,\Pic^0_X) \to H^1(k, \Pic_X)$.
Define $\Br_{1/2}(X)$ to be the subgroup of $\Br_1(X)$ that maps into
the image of $H^1(k,\Pic^0_X)$ in~$H^1(k, \Pic_X)$. (Manin~\cite{Manin}
calls it $\Br'_1(X)$.) In addition, for
$n \ge 1$, let $\Br_{1/2,n}(X)$ be the subgroup of~$\Br_1(X)$ that
maps into the image of $H^1(k, \Pic^0_X)[n]$. Then
\[ \Br_{1/2}(X) = \bigcup_{n \ge 1} \Br_{1/2,n}(X) \]
and
\[ \Ad{X}{k}^{\Br_{1/2}} = \bigcap_{n \ge 1} \Ad{X}{k}^{\Br_{1/2,n}} \,. \]

Recall the definition of $\Ad{X}{k}^{\Alb}$ from Remark~\ref{RemAb} and the
fact that
\[ \Ad{X}{k}^\ab \subset \Ad{X}{k}^{\Alb}
                  = \bigcap_{n \ge 1} \Ad{X}{k}^{\nab} \,. 
\]

\begin{Theorem} \label{BMab}
  Let $X$ be a smooth projective geometrically connected variety,
  and let $n \ge 1$.  Then
  \[ \Ad{X}{k}^{\nab} \subset \Ad{X}{k}^{\Br_{1/2,n}} \,. \]
  In particular,
  \[ \Ad{X}{k}^{\ab} \subset \Ad{X}{k}^{\Alb} \subset \Ad{X}{k}^{\Br_{1/2}} \,.    \]
\end{Theorem}

\begin{Proof}
  Given the first statement, the second statement is clear.
  
  The first statement follows from Thm.~\ref{BMEquality} below. However,
  since our proof of the inclusion given here is fairly simple, we
  include it.
  
  So consider $P \in \Ad{X}{k}^{\nab}$ and $b \in \Br_{1/2,n}(X)$. We have 
  to show that $\langle b, P \rangle_{\Br} = 0$,
  where $\langle \cdot, \cdot \rangle_{\Br}$
  is the Brauer pairing between $\Ad{X}{k}$ and~$\Br(X)$. 
  
  Let $b'$ be the image of~$b$ 
  in~$\Br_1(X)/\Br_0(X) \isom H^1(k, \Pic_X)$, and let 
  $b'' \in H^1(k, \Pic^0_X)[n]$
  be an element mapping to~$b'$ (which exists since $b \in \Br_{1/2,n}(X)$). 
  
  Let $A$ be the Albanese variety of~$X$, and let $V$ be the $k$-torsor 
  under~$A$ that has a canonical map $\phi : X \to V$.
  Then we have $\Pic^0_X \isom \Pic^0_A \isom \Pic^0_V$.
  Since $P \in \Ad{X}{k}^{\nab} \stackrel{\phi}{\to} \Ad{V}{k}^{\nab}$, 
  the latter is nonempty, hence $V$ admits a torsor of the form $(W, A[n])$.
  
  Since $P$ maps into $\Ad{V}{k}^{\nab}$, there is some twist of~$(W, A[n])$
  such that $\phi(P)$ lifts to it. Without loss of generality, 
  $(W, A[n])$ is already this twist, so there is $Q' \in \Ad{W}{k}$ such that 
  $\pi'(Q') = \phi(P)$,
  where $\pi' : W \to V$ is the covering map associated to $(W, A[n])$.
  
  Let $(Y, A[n]) \in \Ab(X)$ be the pull-back of~$(W, A[n])$ to~$X$.
  Then there is some $Q \in \Ad{Y}{k}$ such that $\pi(Q) = P$.
  Now the left hand diagram below induces the one on the right, where
  the rightmost vertical map is multiplication by~$n$:
  \[ \SelectTips{cm}{}
     \xymatrix{ Y \ar[r] \ar[d]_{\pi} & W \ar[d]^{\pi'} & \qquad &
                 \Pic_{Y} &
                 \Pic^0_{Y} \ar[l] & \Pic^0_W \ar[l] \ar@{=}[r] & \Pic^0_A \\
                X \ar[r]^{\phi} & V & \qquad & \Pic_X \ar[u]^{\pi^*} &
                \Pic^0_X \ar[l] \ar[u]^{\pi^*} & 
                     \Pic^0_V \ar[u]^{{\pi'}^*} \ar[l]_{\cong} \ar@{=}[r] & 
                     \Pic^0_A \ar[u]_{\cdot n}
              }
  \]
  Chasing $b''$ around the diagram on the right, after applying $H^1(k, {-})$
  to it, we see that $\pi^*(b') = 0$ in~$\Br(Y)/\Br_0(Y)$. Finally, we have
  \[ \langle b, P \rangle_{\Br} = \langle b', \pi(Q) \rangle_{\Br} 
       = \langle \pi^*(b'), Q \rangle_{\Br} = 0 \,.
  \]
\end{Proof}

So we have the chain of inclusions
\[ \Ad{X}{k}^{\Br} \subset \Ad{X}{k}^{\Br_1} \subset \Ad{X}{k}^{\ab} 
    \subset \Ad{X}{k}^{\Alb} \subset \Ad{X}{k}^{\Br_{1/2}} \,. \]
It is then natural to ask to what extent one might have equality in
this chain of inclusions. We certainly get something when $\Br_{1/2}(X)$
already equals $\Br_1(X)$ or even~$\Br(X)$.

\begin{Corollary} \label{BMabC}
  If $X$ is a smooth projective geometrically connected variety
  such that the canonical map $H^1(k, \Pic^0_X) \to H^1(k, \Pic_X)$
  is surjective, then
  \[ \Ad{X}{k}^{\Br_1} = \Ad{X}{k}^\ab = \Ad{X}{k}^{\Alb} \,. \]
  In particular, if
  $C$ is a curve, then $\Ad{C}{k}^{\Br} = \Ad{C}{k}^\ab$.
\end{Corollary}

\begin{Proof}
  In this case, $\Br_{1/2}(X) = \Br_1(X)$, and so the result follows
  from the two preceding theorems.
  
  When $X = C$ is a curve, then we know that $\Br(C \times_k \bar{k})$
  is trivial (Tsen's Theorem); also $H^1(k, \Pic^0_C) $ surjects onto
  $H^1(k, \Pic_C)$, since the N\'eron-Severi group of~$C$ is~$\Z$ with
  trivial Galois action, and $H^1(k, \Z) = 0$. Hence $\Br(C) = \Br_{1/2}(C)$,
  and the assertion follows.
\end{Proof}

The result of Cor.~\ref{BMabC}
means that we can replace $\Ad{C}{k}^\ab$ by~$\Ad{C}{k}^{\Br}$ everywhere.
For example, from Cor.~\ref{abGen1}, we obtain the following.

\begin{Corollary} \label{Sch}
  Let $C$ be a smooth projective geometrically connected curve over~$k$, and
  let $A$ be its Albanese (or Jacobian) variety.
  Assume that $\Sha(k, A)_{\diw} = 0$.
  \begin{enumerate}\addtolength{\itemsep}{2mm}
    \item If $C$ has a $k$-rational divisor class of degree~$1$ inducing
          a $k$-defined embedding $C \inj A$, then 
          \[ \Ad{C}{k}^{\Br} = \phi^{-1}(\overline{A(k)}) \,, \]
          where $\phi$ denotes 
          the induced map $\Ad{C}{k} \to \Ad{A}{k}$.
    \item If $C$ has no $k$-rational divisor class of degree~$1$,
          then $\Ad{C}{k}^{\Br} = \emptyset$.
  \end{enumerate}
\end{Corollary}

These results can be found in Scharaschkin's thesis~\cite{Scharaschkin}. 
Our approach provides an alternative proof, and the more precise
version in Cor.~\ref{abGen1} shows how to extend the result to the case
when the Shafarevich-Tate group
of the Jacobian is not necessarily assumed to have trivial divisible subgroup.

\medskip

In fact, more is true: we actually have equality in Thm.~\ref{BMab}.

\begin{Theorem} \label{BMEquality}
  Let $X$ be a smooth projective geometrically connected variety. Then
  \[ \Ad{X}{k}^{\nab} = \Ad{X}{k}^{\Br_{1/2,n}} \]
  for all $n \ge 1$. In particular,
  \[ \Ad{X}{k}^{\Alb} = \Ad{X}{k}^{\Br_{1/2}} \,. \]
\end{Theorem}

\begin{Proof}
  This follows from the descent theory of Colliot-Th\'el\`ene and Sansuc.
  Let $M = \Pic^0_X[n]$, and let
  $\lambda : M \to \Pic_X$ be the inclusion. Then the $n$-coverings of~$X$
  are exactly the torsors of type~$\lambda$ in the language of the theory,
  compare for example \cite{Skorobogatov}. (Note that the dual of~$M$
  is $A[n]$, where $A$ is the Albanese variety of~$X$.)
  We have $\Br_\lambda = \Br_{1/2,n}$,
  and the result then follows from Thm.~6.1.2,(a) in~\cite{Skorobogatov}.
\end{Proof}

\begin{Remark}
  Since $\Ad{X}{k}^{\Br_1} \subset \Ad{X}{k}^\ab \subset \Ad{X}{k}^{\Br_{1/2}}$,
  it is natural to ask whether there might be a subgroup $B \subset \Br_1(X)$
  such that $\Ad{X}{k}^\ab = \Ad{X}{k}^B$. As Joost van Hamel pointed out
  to me, a natural candidate for $B$ is the subgroup mapping to the image
  of $H^1(k, \Pic_X^\tau)$ in~$H^1(k, \Pic_X)$, where $\Pic_X^\tau$ is the
  saturation of $\Pic_X^0$ in~$\Pic_X$, i.e., the subgroup of elements
  mapping into the torsion subgroup of the N\'eron-Severi group~$\NS_X$.
  It is tempting to denote this $B$ by $\Br_{2/3}$, but perhaps $\Br_\tau$
  is the better choice. Note that $\Br_\tau = \Br_{1/2}$ when $\NS_X$ is
  torsion free, in which case we have 
  $\Ad{X}{k}^\ab = \Ad{X}{k}^{\Alb} = \Ad{X}{k}^{\Br_{1/2}}$.
\end{Remark}

\begin{Corollary}
  If $C/k$ is a curve that has a rational divisor class of degree~$1$,
  then
  \[ \Ad{C}{k}^{\nab} = \Ad{C}{k}^{\Br[n]} \,. \]
  In words, the information coming from $n$-torsion in the Brauer group
  is exactly the information obtained by an $n$-descent on~$C$.
\end{Corollary}

\begin{Proof}
  Under the given assumptions, 
  $H^1(k, \Pic^0_C) = H^1(k, \Pic_C) = \Br(C)/\Br(k)$,
  and $\Br(k)$ is a direct summand of~$\Br(C)$. Therefore, the images
  of $\Br_{1/2,n}(C)$ and of $\Br(C)[n]$ in $\Br(C)/\Br_0(C)$ agree, 
  and the claim follows.
\end{Proof}

\begin{Corollary}
  If $X$ is a smooth projective geometrically connected variety
  such that the N\'eron-Severi group of~$X$ (over~$\bar{k}$) is
  torsion-free, then there is a finite field extension $K/k$ such that
  \[ \Ad{X}{K}^{\Br_1} = \Ad{X}{K}^\ab \,. \]
\end{Corollary}

\begin{Proof}
  We have an exact sequence
  \[ H^1(k, \Pic^0_X) \To H^1(k, \Pic_X) \To H^1(k, \NS_X) \,. \]
  Since $\NS_X$ is a finitely generated abelian group, 
  the Galois action on it
  factors through a finite quotient $\Gal(K/k)$ of the absolute
  Galois group of~$k$. Then 
  $H^1(K, \NS_X) = \Hom(G_K, \Z^r) = 0$, and the claim
  follows from Thm.~\ref{BMab}.
\end{Proof}

Note that it is not true in general that $\Ad{X}{k}^{\Br_1} = \Ad{X}{k}^\ab$
(even when the N\'eron-Severi group of~$X$ over~$\bar{k}$ is
torsion-free). For example, a smooth cubic surface~$X$ in~$\BP^3$ has 
$\Ad{X}{k}^{\cov} = \Ad{X}{k}$ (since it has trivial geometric fundamental
group), but
may well have $\Ad{X}{k}^{\Br_1} = \emptyset$, even though there are
points everywhere locally. See~\cite{CT-Kanevsky-Sansuc}, where the
algebraic Brauer-Manin obstruction is computed for all smooth diagonal cubic
surfaces 
\[ X : a_1\,x_1^3 + a_2\,x_2^3 + a_3\,x_3^3 + a_4\,x_4^3 = 0 \]
with integral coefficients $0 < a_i < 100$, thereby verifying that it is
the only obstruction against rational points on~$X$ (and thus providing
convincing experimental evidence that this may be true for smooth
cubic surfaces in general).
This computation produces a list of 245 such surfaces with points
everywhere locally, but no rational points, 
since $\Ad{X}{\Q}^{\Br_1} = \emptyset$.

It is perhaps worth mentioning that our condition that $H^1(k, \Pic^0_X)$
surjects onto $H^1(k, \Pic_X)$, which leads to the identification of
the ``algebraic Brauer-Manin obstruction'' and the ``finite abelian descent
obstruction'', is in some sense orthogonal to the situation studied
(quite successfully) in~\cite{CT-Sansuc,CT-Coray-Sansuc,CT-Sansuc-SwD}, 
where it is assumed
that $\Pic_X$ is torsion-free (and therefore $\Pic^0_X$ is trivial), and so 
there can only be a Brauer-Manin obstruction when our condition fails.
There is then no finite abelian descent obstruction, and one has to
look at torsors under tori instead.

\newcommand{\ssu}{\text{\begin{turn}{30}$\subset$\end{turn}}}
\newcommand{\ssd}{\text{\begin{turn}{-30}$\subset$\end{turn}}}
\newcommand{\ssU}{\text{\raisebox{-6pt}{\begin{turn}{30}$\subset$\end{turn}}}}
\newcommand{\ssD}{\text{\raisebox{6pt}{\begin{turn}{-30}$\subset$\end{turn}}}}

In general, we have a diagram of inclusions:
\[ X(k) \subset \overline{X(k)}
   \begin{array}{@{\;}c@{\,}c@{\,}c@{\,}c@{\,}c@{\;}}
     \ssU & \Ad{X}{k}^{\Br} & \subset & \Ad{X}{k}^{\Br_1} & \ssd \\[2mm]
     \ssD & \Ad{X}{k}^{\cov} & \subset & \Ad{X}{k}^{\sol} & \ssu
   \end{array}
   \Ad{X}{k}^\ab \subset \Ad{X}{k}^{\Br_{1/2}}
   \subset \Ad{X}{k}
\]
We expect that every inclusion can be strict. We discuss them in turn.
\begin{enumerate}\addtolength{\itemsep}{2mm}
  \item $X = \BP^1$ has $X(k) \subsetneq \overline{X(k)} = \Ad{X}{k}$.
  \item Skorobogatov's famous example (see \cite{Skorobogatov99} 
        and~\cite{HarariSkorobogatov}) has $\Ad{X}{k}^{\Br} \neq \emptyset$,
        but $\Ad{X}{k}^\sol = \emptyset$, showing that
        $\overline{X(k)} \subsetneq \Ad{X}{k}^{\Br}$ and
        $\Ad{X}{k}^\sol \subsetneq \Ad{X}{k}^\ab$ are both possible.
  \item As mentioned above, \cite{CT-Kanevsky-Sansuc} has examples such
        that $\Ad{X}{k}^{\Br_1} = \emptyset$, but
        $\Ad{X}{k}^{\cov} = \Ad{X}{k}$. This shows that both
        $\overline{X(k)} \subsetneq \Ad{X}{k}^\cov$ and
        $\Ad{X}{k}^{\Br_1} \subsetneq \Ad{X}{k}^\ab$ are possible.
  \item Harari~\cite{Harari1996} has examples, where there is a 
        ``transcendental'', but no ``algebraic'' Brauer-Manin obstruction,
        which means that $\Ad{X}{k}^{\Br} = \emptyset$, but
        $\Ad{X}{k}^{\Br_1} \neq \emptyset$. Hence we can have
        $\Ad{X}{k}^{\Br} \subsetneq \Ad{X}{k}^{\Br_1}$.
  \item If we take a finite nonabelian simple group for $\pi_1(\bar{X})$
        in Cor.~6.1 in~\cite{Harari2000}, then the proof of this result
        shows that $\Ad{X}{k}^\cov \subsetneq \Ad{X}{k}$. On the other
        hand, $\Ad{X}{k}^\sol = \Ad{X}{k}$, since there are only trivial
        torsors in~$\Sol(X)$, compare Lemma~\ref{LemmaPi1}.
  \item It is likely that a construction using Enriques surfaces like that 
        in~\cite{HarariSkorobogatov2} can produce an example such that
        $\Ad{X}{k}^{\Br_{1/2}} = \Ad{X}{k}^{\Alb} = \Ad{X}{k}$,
        since the Albanese variety is trivial, but 
        $\Ad{X}{k}^\ab \subsetneq \Ad{X}{k}$, since there is a nontrivial
        abelian covering.
  \item Finally, in Section~\ref{PropEx} below, we will see many examples 
        of curves~$X$ that have
        $X(k) = \Ad{X}{k}^{\Br_{1/2}} \subsetneq \Ad{X}{k}$.
\end{enumerate}


\bigskip

\subsection*{A new obstruction?}\strut

For curves, we expect the interesting part of the diagram of inclusions
above to collapse: $\overline{X(k)} = \Ad{X}{k}^{\Br_{1/2}}$,
see the discussion in Section~\ref{Conjs} below. For 
higher-dimensional varieties, this is far from true, see the
discussion above. So one could consider a new obstruction obtained
from a combination of the Brauer-Manin and the finite descent obstructions,
as follows. Define
\[ \Ad{X}{k}^{\cov,\Br}
     = \bigcap_{\quad(Y,G) \in \Cov(X)\quad}
         \bigcup_{\xi \in H^1(k,G)} \pi_\xi\Bigl(\Ad{Y_\xi}{k}^{\Br}\Bigr) \,.
\]
(This is similar in spirit to the ``refinement of the Manin obstruction''
introduced in~\cite{Skorobogatov99}.) 

It would be interesting to find out how strong this obstruction is and
whether it is strictly weaker than the obstruction obtained from {\em all}
torsors under (not necessarily finite or abelian) $k$-group schemes. 
Note that the latter is at least as strong as the Brauer-Manin obstruction 
by~\cite[Thm.~4.10]{HarariSkorobogatov} (see also Prop.~5.3.4 
in~\cite{Skorobogatov}), at least if one assumes that all elements
of $\Br(X)$ are represented by Azumaya algebras over~$X$.


\section{Finite descent conditions on curves} \label{PropEx}

Let us now prove some general properties of the notions,
introduced in Section~\ref{CoCoRP} above, of being excellent
w.r.t.\ all, solvable, or abelian coverings in the case of curves. 
In the following, $C$, $D$, etc., will be (smooth projective geometrically 
connected) curves over~$k$. $\iota$ will denote an embedding of~$C$ into its 
Jacobian (if it exists). Also, if $\Ad{C}{k}^{\Br} = \emptyset$
(and therefore $C(k) = \emptyset$, too), we say that {\em the absence
of rational points is explained by the Brauer-Manin obstruction.}
Note that by Cor.~\ref{BMabC}, $\Ad{C}{k}^{\Br} = \Ad{C}{k}^\ab$, which
implies that the absence of rational points is explained by the Brauer-Manin 
obstruction when $C$ is excellent w.r.t.\ abelian coverings and 
$C(k) = \emptyset$. We will use this observation below without explicit
mention.

\begin{Corollary} \label{rk0}
  Let $C/k$ be a curve of genus at least~$1$, with Jacobian~$J$. 
  Assume that $\Sha(k, J)_{\diw} = 0$
  and that $J(k)$ is finite. Then $C$ is excellent w.r.t.\ abelian coverings.
  If $C(k) = \emptyset$, the absence of rational points is explained
  by the Brauer-Manin obstruction.
\end{Corollary}

\begin{Proof}
  By Cor.~\ref{abGen1}, under the assumption on~$\Sha(k, J)$, either
  $\Ad{C}{k}^\ab = \emptyset$, and there is nothing to prove, or else
  \[ \Ad{C}{k}^\ab = \iota^{-1}(\overline{J(k)}) = \iota^{-1}(J(k)) = C(k) \,.
  \]
\end{Proof}

The following result shows that the statement we would like to have 
(namely that $\Ad{C}{k}^\ab = C(k)$) holds for finite subschemes of a curve.

\begin{Theorem} \label{BM0}
  Let $C/k$ be a curve of genus at least~1, and 
  let $Z \subset C$ be a finite subscheme. 
  Then the image of~$\Ad{Z}{k}$ in~$\Ad{C}{k}$ meets $\Ad{C}{k}^{\ab}$ 
  in~$Z(k)$.
  More generally, if $P \in \Ad{C}{k}^\ab$ is such that $P_v \in Z(k_v)$
  for a set of places~$v$ of~$k$ of density~1, then $P \in Z(k)$.
\end{Theorem}

\begin{Proof}
  Let $K/k$ be a finite extension such that $C$ has a rational divisor
  class of degree~1 over~$K$. By Cor.~\ref{abGen1}, we have that
  \[ \Ad{C}{K}^\ab = \iota^{-1}(\Selh(K, J)) \,,\]
  where $\iota : \Ad{C}{K} \to \Ad{J}{K}$ is the map induced by an embedding
  $C \inj J$ over~$K$. Now we apply Thm.~\ref{ZeroDim1} to the image of~$Z$ 
  in~$J$.
  We find that $\iota(P) \in \Selh(K, J)$ and so $\iota(P) \in \iota(Z(K))$.
  Since $\iota$ is injective (even at the infinite places!), we find that
  the image of $P$ in~$\Ad{C}{K}$ is in (the image of) $Z(k)$. Now if $Z(k)$
  is empty, this gives a contradiction and proves the claim in this case.
  Otherwise, $C(k) \supset Z(k)$ is non-empty, and we can take $K = k$ above,
  which gives the statement directly.
\end{Proof}


The following results show that the ``excellence properties'' behave nicely.

\begin{Proposition} \label{Ext}
  Let $K/k$ be a finite extension, and let $C/k$ be a curve of genus
  at least~1. If $C_K$ is excellent w.r.t.\ all, solvable, or abelian
  coverings, then so is~$C$.
\end{Proposition}

\begin{Proof}
  By Prop.~\ref{Extend}, we have
  \[ C(k) \subset \Ad{C}{k}^\cov \subset \Ad{C}{k} \cap \Ad{C}{K}^\cov
          = \Ad{C}{k} \cap C(K) = C(k) \,.
  \]
  Similarly for $\Ad{C}{k}^\sol$ and $\Ad{C}{k}^\ab$.
  Strictly speaking, this means that $C(k)$ and~$\Ad{C}{k}^\cov$ have the
  same image in~$\Ad{C}{K}$. Now, since $C(K)$ has to be finite in order
  to equal~$\Ad{C}{K}^\cov\!$, $C(k)$ is also finite, and we can apply
  Thm.~\ref{BM0} to $Z = C(k) \subset C$ and the set of finite places
  of~$k$. 
\end{Proof}

\begin{Proposition} \label{Cov}
  Let $(D, G) \in \Cov(C)$ (or $\Sol(C)$).
  If all twists~$D_\xi$ of $(D, G)$ are excellent w.r.t.\ all (resp., solvable)
  coverings, then $C$ is excellent w.r.t.\ all (resp., solvable) coverings.
\end{Proposition}

\begin{Proof}
  By~Thm.~\ref{Descent}, $C(k) = \coprod_\xi \pi_\xi(D_\xi(k))$. Now,
  by Prop.~\ref{Cover},
  \[ C(k) \subset \Ad{C}{k}^\cov
          = \bigcup_\xi \pi_\xi\bigl(\Ad{D_\xi}{k}^\cov\bigr)
          = \bigcup_\xi \pi_\xi(D_\xi(k)) = C(k) \,.
  \]
  If $G$ is solvable, the same proof shows the statement for 
  $\Ad{C}{k}^\sol$.
\end{Proof}

\begin{Proposition} \label{Dom}
  Let $C \stackrel{\phi}{\to} X$ be a non-constant morphism over~$k$ from the 
  curve~$C$ into a variety~$X$. If $X$ is excellent w.r.t.\ all, solvable, or
  abelian coverings, then so is~$C$. In particular, if $\Ad{X}{k}^\ab = X(k)$
  and $C(k) = \emptyset$, then the absence of rational points on~$C$ is
  explained by the Brauer-Manin obstruction.
\end{Proposition}

\begin{Proof}
  First assume that $C$ is of genus zero. Then either $\Ad{C}{k} = \emptyset$,
  and there is nothing to prove, or else $C(k)$ is dense in~$\Ad{C}{k}$,
  implying that $X(k) \subsetneq \overline{X(k)} \subset \Ad{X}{k}^\cov$ and 
  thus contradicting the assumption.
  
  Now assume that $C$ is of genus at least~1.
  Let $P \in \Ad{C}{k}^{\cov/\sol/\ab}$. Then by Thm.~\ref{MorIncl}, 
  $\phi(P) \in \Ad{X}{k}^{\cov/\sol/\ab} = X(k)$. Let $Z \subset C$ be the 
  preimage (subscheme) of $\phi(P) \in X(k)$ in~$C$. This is finite,
  since $\phi$ is non-constant. Then we have
  that $P$ is in the image of~$\Ad{Z}{k}$ in~$\Ad{C}{k}$. Now apply
  Thm.~\ref{BM0} to conclude that
  \[ P \in \Ad{C}{k}^\ab \cap \Ad{Z}{k} = Z(k) \subset C(k) \,. \]
\end{Proof}


As an application, we have the following.

\begin{Theorem} \label{CorMor}
  Let $C \to A$ be a non-constant morphism over~$k$ of a curve~$C$ into
  an abelian variety~$A$. Assume that $\Sha(k, A)_{\diw} = 0$ and that 
  $A(k)$ is finite. (For example, this is the case when $k = \Q$ and
  $A$ is modular of analytic rank zero.)
  Then $C$ is excellent w.r.t.\ abelian coverings.
  In particular, if $C(k) = \emptyset$, then the absence of rational points
  on~$C$ is explained by the Brauer-Manin obstruction.
\end{Theorem}

\begin{Proof}
  By Cor.~\ref{AV}, we have $\Ad{A}{k}^\ab = A(k)$. Now by Prop.~\ref{Dom}, 
  the claim follows.
\end{Proof}

This generalizes a result proved by Siksek~\cite{Siksek} under
additional assumptions on the Galois action on the fibers of~$\phi$
above $k$-rational points of~$A$, in the case that $C(k)$ is empty. 
A similar observation was made 
independently by Colliot-Th\'el\`ene~\cite{CTremarque}. Note that both
previous results are in the context of the Brauer-Manin obstruction.

\begin{Examples}
  We can use Thm.~\ref{CorMor} to produce many examples of curves $C$
  over~$\Q$ that are excellent w.r.t.\ abelian coverings. Concretely, let 
  us look at the curves $C_a : y^2 = x^6 + a$, where $a$ is a non-zero
  integer. $C_a$ maps to the two elliptic curves $E_a : y^2 = x^3 + a$
  and $E_{a^2}$ (the latter by sending $(x,y)$ to $(a/x^2, ay/x^3)$).
  So whenever one of these elliptic curves has (analytic) rank zero,
  we know that $C_a$ is excellent w.r.t.\ abelian coverings. 
  For example, this is the
  case for all $a$ such that $|a| \le 20$, with the exception of
  $a = -15, -13, -11, 3, 10, 11, 15, 17$. Note that 
  $C_a(\Q)$ is always non-empty (there are two rational points at infinity).
\end{Examples}

We can even show a whole class of interesting curves to be excellent
w.r.t.\ abelian coverings.

\begin{Corollary} \label{ModCurves}
  If $C/\Q$ is one of the modular curves $X_0(N)$, $X_1(N)$, $X(N)$ and
  such that the genus of~$C$ is positive, then $C$ is excellent 
  w.r.t.\ abelian coverings.
\end{Corollary}

\begin{Proof}
  By a result of Mazur~\cite{Mazur}, every Jacobian $J_0(p)$ of~$X_0(p)$,
  where $p = 11$ or $p \ge 17$ is prime, has a nontrivial factor of
  analytic rank zero. Also, if $M \mid N$, then there are nonconstant
  morphisms $X_1(N) \to X_0(N) \to X_0(M)$, hence the assertion is true
  for all $X_0(N)$ and $X_1(N)$ such that $N$ is divisible by one of the
  primes in Mazur's result. For the other minimal $N$ such that $X_0(N)$
  (resp., $X_1(N)$) is of positive genus, William Stein's 
  tables~\cite{SteinTables} prove that there is a factor of $J_0(N)$
  (resp., $J_1(N)$) of analytic rank zero. So we get the result for
  all $X_0(N)$ and $X_1(N)$ of positive genus. Finally, $X(N)$ maps to
  $X_0(N^2)$, and so we obtain the result also for $X(N)$ (except in
  the genus zero cases $N = 1, 2, 3, 4, 5$).
\end{Proof}

For another example, involving high-genus Shimura curves, 
see~\cite{SkorobogatovSh}.

\begin{Remark} \label{RkSection}
  There is some relation with the ``Section Conjecture'' from 
  Grothen\-dieck's anabelian geometry~\cite{Grothendieck}. 
  Let $C/k$ be a smooth projective geometrically connected
  curve of genus~$\ge 2$.
  One can prove that if $C$ has the ``section property'', then $C$ is
  excellent w.r.t.\ all coverings, which in turn implies that $C$ has
  the ``birational section property''.
  See Koenigsmann's paper~\cite{Koenigsmann} for definitions.
  For example, all the curves $X_0(N)$, $X_1(N)$ and $X(N)$ have the
  birational section property if they are of higher genus.
\end{Remark}


\section{Discussion} \label{Conjs}

In the preceding section, we have seen that we can construct many 
examples of higher-genus curves that are excellent w.r.t.\ abelian
coverings. This leads us to state the following conjecture.

\begin{Conjecture}[Main Conjecture] \label{MainConj}
  If $C$ is a smooth projective geometrically connected curve over a 
  number field~$k$, then $C$ is very good.
\end{Conjecture}

By what we have seen, for curves of genus~$1$, this is equivalent to
saying that the divisible subgroup of $\Sha(k, E)$ is trivial, for
every elliptic curve $E$ over~$k$. For curves $C$ of higher genus,
the statement is equivalent to saying that $C$ is excellent w.r.t.\ abelian coverings.
We recall that our conjecture would follow in this case from the
``Adelic Mordell-Lang Conjecture'' formulated in Question~\ref{AML}.

\begin{Remark}
  When $k$ is a global function field of characteristic~$p$, then the Main 
  Conjecture holds
  when $J = \Jac_C$ has no isotrivial factor and $J(k^{\text {sep}})[p^\infty]$
  is finite. 
  See recent work by Poonen and Voloch~\cite{PoonenVoloch}.
\end{Remark}

If the Main Conjecture holds for~$C$ and $C(k)$ is empty,
then (as previously discussed) we can find a torsor that has no twists with 
points everywhere locally and thus {\em prove} that $C(k)$ is empty.
The validity of the conjecture (even just in case $C(k)$ is
empty) therefore implies that {\em we can algorithmically decide whether
a given smooth projective geometrically connected curve over a 
number field~$k$ has rational points or not.}

In Section~\ref{BM} above, we have shown that for a curve~$C$, we have
\[ \Ad{C}{k}^{\ab} = \Ad{C}{k}^{\Br} \,, \]
where on the right hand side, we have the {\em Brauer subset} of~$\Ad{C}{k}$,
i.e., the subset cut out by conditions coming from the Brauer group of~$C$.
One says that there is a {\em Brauer-Manin obstruction} against rational
points on~$C$ if $\Ad{C}{k}^{\Br} = \emptyset$. A corollary 
of our Main Conjecture is that the Brauer-Manin obstruction is the
only obstruction against rational points on curves over number fields
(which means that $C(k) = \emptyset$ implies $\Ad{C}{k}^{\Br} = \emptyset$).
To our knowledge, before this work (and Poonen's heuristic, see his conjecture
below, which was influenced by discussions we had at the IHP in Paris in 
Fall~2004) nobody gave a conjecturally positive answer to the question,
first formulated on page~133 in~\cite{Skorobogatov}, whether the Brauer-Manin
obstruction might be the only obstruction against rational points on curves.
No likely counter-example is known, but
there is an ever-growing list of examples, for which
the failure of the Hasse Principle could be explained by the Brauer-Manin
obstruction; see the discussion below (which does not pretend to be
exhaustive) or also Skorobogatov's
recent paper~\cite{SkorobogatovSh} on Shimura curves.

\medskip

Let $v$ be a place of~$k$. Under a {\em local condition at~$v$} on
a rational point $P \in C$, we understand the requirement that the
image of $P$ in~$C(k_v)$ is contained in a specified closed and open
(``clopen'') subset of~$C(k_v)$. If $v$ is an infinite place, this
just means that we require $P$ to be on some specified connected component(s)
of~$C(k_v)$; for finite places, we can take something like a ``residue class''.
With this notion, the Main Conjecture~\ref{MainConj} above is equivalent
to the following statement.

{\em Let $C/k$ be a curve as above. Specify local conditions at finitely
many places of~$k$ and assume that there is no point in~$C(k)$ satisfying
these conditions. Then there is some $n \ge 1$ such that no point 
in $\prod_v X_v \subset \Ad{C}{k}$ survives the $n$-covering of~$C$,
where $X_v$ is the set specified by the local condition at those
places where a condition is specified, and $X_v = C(k_v)$ (or $\pi_0(C(k_v))$)
otherwise.}

This says that the ``finite abelian'' obstruction (equivalently, the
Brauer-Manin obstruction) is the only obstruction against weak approximation
in~$\Ad{C}{k}$.

We see that the conjecture implies that we can decide if a given finite
collection of local conditions can be satisfied by a rational point.
Now the question is how practical it might be to actually do this in 
concrete cases. For certain classes of curves and specific values of~$n$,
it may be possible to explicitly and efficiently find the relevant twists.
For example, this can be done for hyperelliptic curves and $n = 2$, 
compare~\cite{BruinStoll2D}. However, for general curves and/or general~$n$,
this approach is likely to be infeasible.

On the other hand, assume that we can find $J(k)$ 
explicitly, where $J$, as usual, is the Jacobian of~$C$. This is the
case (at least in principle) when $\Sha(k, J)_{\diw} = 0$. Then we
can approximate $\Ad{C}{k}^\ab$ more and more precisely by looking
at the images of~$\Ad{C}{k}$ 
and of~$J(k)$ in $\prod_{v \in S} J(k_v)/N J(k_v)$
for increasing~$N$ and finite sets~$S$ of places of~$k$. If $C(k)$
is empty and the Main Conjecture holds, then for some choice of~$S$ and~$N$,
the two images will not intersect, giving an explicit proof that
$C(k) = \emptyset$. An approach like this was proposed (and carried out
for some twists of the Fermat quartic) by
Scharaschkin~\cite{Scharaschkin}. See~\cite{Flynn} for an implementation of 
this method and~\cite{BruinStollBM} for improvements.
In~\cite{PSS}, this procedure is used to rule out rational points
satisfying certain local conditions on a genus~$3$ curve whose Jacobian
has Mordell-Weil rank~$3$.

In order to test the conjecture, Nils Bruin and the author conducted 
an experiment, see~\cite{BruinStollEx}. We considered all genus~$2$ curves
over~$\Q$ of the form
\begin{equation} \label{SmallG2}
  y^2 = f_6\,x^6 + f_5\,x^5 + \dots + f_1\,x + f_0 
\end{equation}
with coefficients $f_0, \dots, f_6 \in \{-3, -2, \dots, 3\}$. For each
isomorphism class of curves thus obtained, we attempted to decide if
there are rational points or not. On about 140\,000 of these roughly 200\,000
curves (up to isomorphism), we found a (fairly) small rational point.
Of the remaining about 60\,000, about half failed to have local points at some
place. On the remaining about 30\,000 curves, we performed a $2$-descent
and found that for all but 1\,492 curves~$C$, 
$\Ad{C}{\Q}^{\text{$2$-ab}} = \emptyset$, proving that $C(\Q) = \emptyset$
as well. For the 1\,492 curves that were left over, we found generators
of the Mordell-Weil group (assuming the Birch and Swinnerton-Dyer Conjecture
for a small number of them) and then did a computation along the lines
sketched above. This turned out to be successful for {\em all} curves,
proving that none of them has a rational point. The conclusion is that
the Main Conjecture holds for curves $C$ as in~\eqref{SmallG2} if
$C(\Q) = \emptyset$, assuming $\Sha(\Q, J)_{\diw} = 0$ for the Jacobian~$J$
if $C$ is one of the 1\,492 curves mentioned, and assuming in addition 
the Birch and Swinnerton-Dyer Conjecture if $C$ is one of 42 specific curves.

\medskip

At least in case $C(k)$ is empty, there are heuristic arguments 
due to Poonen~\cite{PoonenHeur} that
suggest that an even stronger form of our conjecture might be true.

\begin{Conjecture}[Poonen] \label{PoonenConj}
  Let $C$ be a smooth projective geometrically connected curve of genus 
  $\ge 2$ over a number field~$k$, and assume that $C(k) = \emptyset$.
  Assume further that $C$ has a rational divisor class of degree~$1$,
  and let $\iota : C \to J$ be the induced embedding of $C$ into its
  Jacobian~$J$.
  Then there is a finite set~$S$ of finite places of good reduction for~$C$
  such that the image of $J(k)$ in $\prod_{v \in S} J(\F_v)$
  does not meet $\prod_{v \in S} \iota(C(\F_v))$.
\end{Conjecture}

Note that under the assumption $\Sha(k, J)_{\diw} = 0$, we must have
a rational divisor (class) of degree~1 on~$C$ whenever 
$\Ad{C}{k}^\ab \neq \emptyset$, compare Cor.~\ref{abGen1}, so the
condition above is not an essential restriction.

Let us for a moment assume that Poonen's Conjecture holds and that
all abelian varieties $A/k$ satisfy $\Sha(k, A)_{\diw} = 0$. Then for all
curves~$C/k$ of higher genus, $C(k) = \emptyset$ implies 
$\Ad{C}{k}^{\ab} = \emptyset$. If we apply this observation to
coverings of~$C$, then we find that $C$ must be excellent w.r.t.\ solvable
coverings. The argument goes like this. Let $P \in \Ad{C}{k}^{\sol}$,
and assume $P \notin C(k)$. There are only finitely many rational points
on~$C$, hence there is an~$n$ such that $P$ lifts to a
different $n$-covering $D$ of~$C$ than all the rational points. 
(Take $n$ such that $P - Q$ is not divisible by~$n$ in $\Ad{J}{k}$,
for all $Q \in C(k)$, where $J$ is the Jacobian of~$C$.)
In particular, $D(k)$ must be empty. But then, by Poonen's Conjecture,
we have $\Ad{D}{k}^{\ab} = \emptyset$, so $P$ cannot lift to~$D$ either.
This contradiction shows that $P$ must be a rational point.

In particular, this would imply that all higher-genus curves have the
`birational section property', compare Remark~\ref{RkSection}.

\medskip

A more extensive and detailed discussion of these conjectures, their
relations to other conjectures, and evidence for them will be published 
elsewhere.


\end{document}